\documentclass[12pt]{amsart}
\usepackage{amssymb}

\begin{document}

\def\1#1{\overline{#1}}
\def\2#1{\widetilde{#1}}
\def\3#1{\widehat{#1}}
\def\4#1{\mathbb{#1}}
\def\5#1{\frak{#1}}
\def\6#1{{\mathcal{#1}}}

\def\C{{\4C}}
\def\bC{{\4C}}
\def\R{{\4R}}
\def\bR{{\4R}}
\def\M{{\3M}}
\def\N{{\4N}}
\def\T{{\sf T}}
\def\bZ{{\4Z}}


\def\cn{{\C^n}}
\def\cnn{{\C^{n'}}}
\def\ocn{\2{\C^n}}
\def\ocnn{\2{\C^{n'}}}
\def\ol#1\overline{#1}


\def\const{{\rm const}}
\def\rk{{\rm rank\,}}
\def\id{{\sf id}}
\def\aut{{\sf aut}}
\def\Aut{{\sf Aut}}
\def\CR{{\rm CR}}
\def\GL{{\sf GL}}
\def\J{{\sf J}}
\def\Re{{\sf Re}\,}
\def\re{{\sf Re}\,}
\def\Im{{\sf im}\,}
\def\im{{\sf Im}\,}
\def\ker{{\sf ker}\,}

\def\codim{{\rm codim}}
\def\crd{\dim_{{\rm CR}}}
\def\crc{{\rm codim_{CR}}}

\def\phi{\varphi}
\def\eps{\varepsilon}
\def\d{\partial}
\def\a{\alpha}
\def\b{\beta}
\def\g{\gamma}
\def\G{\Gamma}
\def\D{\Delta}
\def\Om{\Omega}
\def\k{\kappa}
\def\l{\lambda}
\def\L{\Lambda}
\def\z{{\bar z}}
\def\w{{\bar w}}
\def\Z{{\1Z}}
\def\t{\theta}


\newtheorem{Thm}{Theorem}[section]
\newtheorem{Cor}[Thm]{Corollary}
\newtheorem{Pro}[Thm]{Proposition}
\newtheorem{Lem}[Thm]{Lemma}
\theoremstyle{definition}\newtheorem{Def}[Thm]{Definition}
\newtheorem{Rem}[Thm]{Remark}
\newtheorem{Exa}[Thm]{Example}
\newtheorem{Exs}[Thm]{Examples}

\def\bl{\begin{Lem}}
\def\el{\end{Lem}}
\def\bp{\begin{Prop}}
\def\ep{\end{Prop}}
\def\bt{\begin{Thm}}
\def\et{\end{Thm}}
\def\bc{\begin{Cor}}
\def\ec{\end{Cor}}
\def\bd{\begin{Def}}
\def\ed{\end{Def}}
\def\br{\begin{Rem}}
\def\er{\end{Rem}}
\def\be{\begin{Exa}}
\def\ee{\end{Exa}}
\def\bpf{\begin{proof}}
\def\epf{\end{proof}}
\def\ben{\begin{enumerate}}
\def\een{\end{enumerate}}

\def\cn{{\C^N}}
\def\cnn{{\C^{N'}}}

\title[Jet parametrization of local CR automorphisms]
{Finite jet determination of local analytic CR automorphisms
and their parametrization by $2$-jets in the finite type case}
\author[P. Ebenfelt, B. Lamel, and  D. Zaitsev]{Peter Ebenfelt,
Bernhard Lamel, and Dmitri Zaitsev}
\footnotetext{{\rm The first
  author is a Royal Swedish Academy of Sciences Research Fellow
  supported by a grant from the Knut and Alice Wallenberg
Foundation.\newline}} \footnotetext{{\rm The second author is
supported by the ANACOGA research network.\newline}}
\address{P. Ebenfelt: Department of Mathematics, University
of California at San Diego, La Jolla, CA 92093-0112, USA}
\email{pebenfel@math.ucsd.edu }
\address{B. Lamel: Department of Mathematics, University
of Illinois at Urbana-Champaign, Urbana, IL 61801, USA}
\email{blamel@math.ucsd.edu }
\address{D. Zaitsev: Dipartimento di Matematica,
Universit\`a di Padova, via G. Belzoni 7, 35131 Padova, ITALY}
\email{zaitsev@math.unipd.it}

\abstract We show that germs of local real-analytic CR
automorphisms of a real-analytic hypersurface $M$ in $\C^2$ at a
point $p\in M$ are uniquely determined by their jets of some
finite order at $p$ if and only if $M$ is not Levi-flat near $p$.
This seems to be the first necessary and sufficient result on
finite jet determination and the first result of this kind in the
infinite type case.

If $M$ is of finite type at $p$, we prove a stronger assertion:
the local real-analytic CR automorphisms of $M$ fixing $p$ are
analytically parametrized (and hence uniquely determined) by their
$2$-jets at $p$. This result is optimal since the automorphisms of
the unit sphere are not determined by their $1$-jets at a point of
the sphere.

We also give an application to the dynamics of germs of local
biholomorphisms of $\C^2$.
\endabstract

\def\Label#1{\label{#1}}

\maketitle

\section{Introduction}\Label{intro}
By {\sc H. Cartan}'s classical uniqueness theorem \cite{CaH35}, a
biholomorphic automorphism of a bounded domain $D\subset\C^n$ is
uniquely determined by its value and its first order derivatives
(that is, by its $1$-jet) at any given point $p\in D$. The example
of the unit ball $D$ shows that if $p$ is taken on the boundary
$\d D$ the same uniqueness phenomenon does not hold for $1$-jets
but rather for $2$-jets at $p$. (In this case, all automorphisms
extend holomorphically across the boundary and, hence, any jet at
a boundary point is well defined.). More generally, results of
{\sc E.~Cartan} \cite{Ca32a, Ca32b}, {\sc N.~Tanaka} \cite{Ta62}
and {\sc S.-S. Chern - J.K. Moser} \cite{CM} (see also {\sc
H.~Jacobowitz} \cite{J77}) show that unique determination by
$2$-jets at a point $p$ holds for (germs at $p$ of) local
biholomorphisms of $\C^n$ sending a (germ at $p$ of an) open piece
$M\subset\d D$ into itself provided $M$ is a {\em
Levi-nondegenerate} real-analytic hypersurface.

The case of a {\em degenerate} real-analytic hypersurface $M$ is
much less understood, even in $\C^2$. It was previously known that
unique determination, as above, by jets at $p$ of some finite
order $k$ holds if $M$ is of {\em finite type} at $p$, due to a
recent result of the first author jointly with {\sc M.S.~Baouendi}
and {\sc L.P.~Rothschild} \cite[Corollary~2.7]{JAMS}. (We remark
here that both notions of finite type, i.e.\ that in the sense of
\cite{Ko72,BG77} and that of \cite{D82},
coincide in $\bC^2$.) On the other hand,
the only known situation where the germs of local biholomorphisms
of $\C^2$ sending $M$ into itself are not uniquely determined by
their $k$-jets at $p$, for any $k$, is that where $M$ is {\em Levi
flat}. The first main result of this paper fills the gap between
these two situations by showing that the Levi-flat case is indeed
the only exception:

\begin{Thm}\Label{main}
Let $M\subset\C^2$ be a real-analytic hypersurface which is not
Levi-flat near a point $p\in M$. Then there exists an integer $k$
such that if $H^1$ and $H^2$ are germs at $p$ of local
biholomorphisms of $\C^2$ sending $M$ into itself with $j^k_p
H^1=j^k_p H^2$, then $H^1\equiv H^2$.
\end{Thm}

We remark that this and all the following results in this paper
about local biholomorphisms sending $M$ into itself also hold for
local biholomorphisms sending $M$ into another real-analytic
hypersurface $M'\subset\C^2$. Indeed, any fixed local
biholomorphism $f_0$ sending $M$ into $M'$ defines a one-to-one
correspondence between the set of germs preserving $M$ and that
sending $M$ to $M'$ (as well as those of their jets) via $g\mapsto
f_0\circ g$.

Since germs of local biholomorphisms of $\C^2$ sending $M$ into
itself are in one-to-one correspondence with germs of
real-analytic CR automorphisms of $M$ by a theorem of {\sc
Tomassini} \cite{To66}, we obtain the following immediate
consequence of Theorem~\ref{main}:

\bc Let $M\subset\C^2$ be a real-analytic hypersurface which is
not Levi-flat near a point $p\in M$. Then there exists an integer
$k$ such that if $h^1$ and $h^2$ are germs at $p$ of real-analytic
CR automorphisms of $M$ with $j^k_p h^1=j^k_p h^2$, then
$h^1\equiv h^2$. \ec

The proof of Theorem \ref{main} relies on
the parametrization of local biholomorphisms
along the zero Segre variety given in \S\ref{sec-jet}, the method
of singular complete systems recently developed by the first author
\cite{Enonm}
and on a finite determination result for solutions of singular
differential equations given in \S\ref{sec-ode}.

The second main result of this paper improves the results in the
{\em finite type case} mentioned above in two different
directions. First we show that, just as in the Levi-nondegenerate
case, the $2$-jets are always sufficient for unique determination.
This conclusion contrasts strongly with most known results for
hypersurfaces of finite type, where one usually has to take at
least as many derivatives as the type (i.e.\ the minimal length of
commutators of vector fields on $M$ in the complex tangent
direction required to span the full tangent space). For instance,
in the setting of Theorem~\ref{main}, upper estimates on how large
the number $k$ must be chosen were previously only known for {\em
finitely nondegenerate} manifolds $M$, due to results in
\cite{Asian} (see also \cite{Z1} and \cite{CAG}) and the estimate
for $k$ in this case was at least twice the type of $M$ minus two.
Although we shall not define finite nondegeneracy in this paper
(but refer the reader instead e.g.\ to the book \cite{BER}), we
would like to point out that finite nondegeneracy is a strictly
stronger notion than finite type; e.g.\ the hypersurface
$$M:=\{(z,w)\in\C^2 : \im w = |z|^4\}$$ is of finite type but is
not finitely nondegenerate at $0$. In this paper, we show that
$2$-jets are sufficient for unique determination as in
Theorem~\ref{main} regardless of the type of $M$ at $p$.

The second improvement in the finite type case is the stronger
conclusion (than that of Theorem~\ref{main}) that local
biholomorphisms are not only uniquely determined but are
analytically parametrized by their $2$-jets. We denote by
$G^2_p(\C^2)$ the group of all $2$-jets at $p$ of local
biholomorphisms $H\colon(\C^2,p)\to(\C^2,p')$.

\begin{Thm}\Label{main-jetpar}
Let $M\subset\bC^2$ be a real-analytic hypersurface of finite
type at a point $p\in M$.
Then there exist an open subset $\Omega\subset\C^2\times G^2_{p}(\C^2)$
and a real-analytic map $\Psi(Z,\L)\colon\Omega\to \C^2$, which in
addition is holomorphic in $Z$, such that the following holds.
For every local biholomorphism $H$ of $\C^2$ sending $(M,p)$ into
itself,
the point $(p,j^2_p H)$ belongs to $\Omega$ and the identity
\begin{equation*}
H(Z) \equiv \Psi (Z, j_p^2 H)
\end{equation*}
holds for all $Z\in\C^2$ near $p$.
\end{Thm}

The conclusion of Theorem~\ref{main-jetpar} was previously known
for Levi-nondegenerate hypersurfaces in $\C^N$ due to the
classical results mentioned above and for Levi-nondegenerate CR
submanifolds of higher codimension due to a more result of {\sc
V.K. Beloshapka} \cite{Be90}. The existence of a $k$-jet
parametrization, for some $k$, was known in the more general case
(than that of Levi nondegenerate hypersurfaces) where $M$ is
finitely nondegenerate at $p$ (see \cite{MA} and previous results
in \cite{Asian} and \cite{Z1}). For other results on finite jet
determination and finite jet parametrization of local CR
automorphisms (in the finite type case), the reader is referred to
the papers \cite{Han2, CAG, Hay, MA, Lamel, BMR, Ejet, Ki01, KZ}.
Theorem~\ref{main-jetpar} will be a consequence of
Theorem~\ref{thm-jetpar} which will be proved in
\S\ref{sec-finite}.

We conclude this introduction by giving some applications of our
main results. In view of a regularity result of {\sc X. Huang}
\cite{HCPDE}, a unique determination result can be formulated for
{\em continuous} CR homeomorphisms as follows. Recall that a
continuous mapping $h\colon M\to \bC^2$ is called CR if it is
annihilated, in the sense of distributions, by all CR vector
fields on $M$ (i.e.\ by the traces of $(0,1)$-vector fields in
$\C^2$ tangent to $M$). A homeomorphism $h$ between two
hypersurfaces $M$ and $M'$ is called CR if both $h$ and $h^{-1}$
are CR mappings. {\sc X.~Huang} showed in the above mentioned
paper that every continuous CR mapping between two real-analytic
hypersurfaces $M, M'\subset \C^2$ of finite type extends
holomorphically to a neighborhood of $M$ in $\C^2$. Thus
Theorem~\ref{main-jetpar} implies:

\begin{Cor}\Label{rem-finite1}
Let $M,M'\subset\bC^2$ be real-analytic
 hypersurfaces of finite type.
Then, for any $p\in M$, if $h_1$ and $h_2$
are germs at $p$ of local CR homeomorphisms between $M$ and $M'$ such
that
\begin{equation}
h_1(x)-h_2(x)=o(\|x\|^2), \quad x\to 0,
\end{equation}
then $h_1\equiv h_2$,
where $x=(x_1,x_2,x_3)$ are any local coordinates on $M$ vanishing at
$p$.
\end{Cor}

Observe that, since a homeomorphism need not preserve the
vanishing order, a uniqueness statement in the spirit of
Corollary~\ref{rem-finite1} may not be reduced to the case $M=M'$
in general. In our case, however, the reduction is possible due to
the above mentioned theorem of {\sc X.~Huang}.

Our next application of Theorem~\ref{main-jetpar} is a structure
result for the group $\Aut(M,p)$ of all local biholomorphisms
$H\colon (\C^2,p)\to (\C^2,p)$ sending $M$ (with $p\in M$) into
itself. A fundamental problem, usually referred to as the {\em
local biholomorphic equivalence problem}, is to determine for
which pairs of (germs of) real submanifolds $(M,p)$ and $(M',p')$
there exist local biholomorphisms sending $(M,p)$ into $(M',p')$
or, formulated in a slightly different way, to describe, for a
given germ of a manifold $(M,p)$, its equivalence class under
local biholomorphic transformations. There is, of course, no loss
of generality in assuming that $p=p'$.
The action of the group $\6E_{p}$ of all local biholomorphisms
$H\colon (\bC^2,p)\to (\bC^2,p)$ will take us from any germ
$(M,p)$ to any other germ $(M',p)$ in the same equivalence class.
However, $\mathcal E_{p}$ does not in general act freely on the
equivalence class of $(M,p)$. To understand the structure of the
equivalence classes one is therefore led to study the structure of
the isotropy group $\Aut(M,p)$. Observe that $\Aut(M,p)$ is a
topological group equipped with a natural direct limit topology
which it inherits as a subgroup of $\6E_{p}$. A sequence of germs
$H^j$ is convergent if all germs $H^j$ extend holomorphically to a
common neighborhood of $p$ on which they converge uniformly (c.f.\
e.g. \cite{Asian}). By standard techniques (see e.g.\ \cite{Asian}
and \cite{MA}), Theorem~\ref{main-jetpar} implies the following:

\bc\Label{thm-finite2}
Let $(M,p)$ be a germ of a real-analytic hypersurface in $\bC^2$ of
finite type.
Then the jet evaluation homomorphism
\begin{equation}\Label{eq-finite2}
j^2_{p}\colon \Aut(M,p)\to G_{p}^2(\C^2)
\end{equation}
is a  homeomorphism onto a closed Lie subgroup of $G_{p}^2(\C^2)$
and hence defines a Lie group structure on $\Aut(M,p)$. \ec

We conclude with an application to the dynamics of germs of local
biholomorphisms of $\C^2$:

\bt\Label{dynamics} Let $H\colon (\C^2,0)\to(\C^2,0)$ be a local
biholomorphism tangent to the identity at $0$, i.e.\ of the form
$H(Z)=Z+O(|Z|^2)$. Suppose that $H$ preserves a germ of a
Levi-nonflat real-analytic hypersurface at $0$. Then $H$ fixes
each point of a complex hypersurface through $0$. \et

Theorem~\ref{dynamics} is a consequence of
Theorem~\ref{thm-jetsegre}, whose statement and proof are given in
\S\ref{sec-jet}.

\section{Preliminaries}\Label{sec-preli}

Recall that a real-analytic hypersurface $M\subset\C^N$ ($N\ge 2$)
is called {\em Levi-flat} if its Levi form
$$L(\xi):=\sum_{k,j} \frac{\d^2\rho}{\d z_k \d \bar z_j} \xi_k
\bar\xi_j$$
vanishes identically on the complex tangent subspace
$T^c_p M:= T_p M\cap iT_p M$ for any $p\in M$,
where $M$ is locally given by $\{\rho = 0\}$ with $d\rho\ne 0$.
It is well-known (and not difficult to see)
that $M$ is Levi-flat if and only if, at any point $p\in M$, there are
local holomorphic coordinates
$(z,w)\in\C^{N-1}\times\C$ in which $M$ has the form $\{\im w=0\}$.

In general, let $M\subset\bC^N$ be a real-analytic hypersurface
with $p\in M$. We may choose local coordinates $(z,w)\in
\bC^{N-1}\times \bC$, vanishing at $p$, so that $M$ is defined
locally near $p=(0,0)$ by an equation of the form
\begin{equation}\Label{eq-normal1}
M :\quad\im w=\phi(z,\bar z,\re w),
\end{equation}
where $\phi(z,\bar z, s)$ is
a real-valued, real-analytic function satisfying
\begin{equation}\Label{id-normal-phi}
\phi(z,0,s)\equiv\phi(0,\chi,s)\equiv 0.
\end{equation}
Such coordinates are called {\it normal
coordinates} for $M$ at $p$; the reader is referred to e.g.\
\cite{BER} for the existence of
such coordinates and related basic material concerning real
submanifolds in
complex space. We mention also that the hypersurface $M$ is of finite
type at $p=(0,0)$ if and only if, in normal coordinates,
$\phi(z,\chi,0)\not\equiv 0$.

\section{Parametrization of jets along the zero Segre
variety}\Label{sec-jet}

Let $M\subset\bC^N$ be a real-analytic hypersurface with $p\in M$.
We shall choose normal coordinates $(z,w)\in \bC^{N-1}\times \bC$
for $M$ at $p$; i.e.\ $(z,w)$ vanishes at $p$ and $M$ is defined
locally near $p=(0,0)$ by (\ref{eq-normal1}) where $\phi(z,\bar z,
s)$ is a real-valued, real-analytic function satisfying
(\ref{id-normal-phi}).  Denote by $J^k_{0,0}(\bC^N)$ the space of
all $k$-jets at $0$ of holomorphic mappings $H\colon(\bC^N,0)\to
(\bC^N,0)$, and by $\T(\bC^N)\subset J^1_{0,0}(\bC^N)$ the group
of invertible upper triangular matrices. Given coordinates $Z$ and
$Z'$ near the origins of the source and target copy of $\bC^N$,
respectively, we obtain associated coordinates
$\Lambda=(\Lambda_i^{\alpha})_{1\leq i\leq N, 1\leq |\alpha|\leq
k}\in J^k_{0,0}(\bC^N)$, (where $i\in \mathbb Z_+$ and
$\alpha\in\mathbb Z_+^N$), in which a jet $j^k_0H$ is given by
$\Lambda_i^{\alpha}=\partial_Z^\alpha H_i(0)$. In this paper, we
shall be concerned with the situation $N=2$ where we have
coordinates $(z,w)$ near $0$ in the source copy of $\bC^2$ and
$(z',w')$ near 0 in the target $\bC^2$. A map $H\colon
(\bC^2,0)\to (\bC^2,0)$ is then given in coordinates by
$H(z,w)=(F(z,w), G(z,w))$. We shall use, for a given $k$, the
notation $\Lambda=(\lambda^{ij},\mu^{ij})_{1\leq i+j\leq k}$, for
the associated coordinates on $J^k_{0,0}(\bC^2)$, where
$\lambda^{ij}=F_{z^i w^j}(0)$ and $\mu^{ij}=G_{z^i w^j}(0)$, and
we use the notation $F_{z^i w^j}=\d^i_z\d^j_w F$ etc.\ for partial
derivatives. In this section, we shall prove the following result.

\begin{Thm}\Label{thm-jetsegre}
Let $M,M'\subset\bC^2$ be real-analytic hypersurfaces that are not
Levi-flat, and let $(z,w)\in\bC^2$ and $(z',w')\in \bC^2$ be
normal coordinates for $M$ and $M'$ vanishing at $p\in M$ and
$p'\in M'$, respectively. Then, for any integer $k\ge 0$, the
identity
\begin{equation}\Label{eq-jetsegre0}
H_{w^k}(z,0) \equiv \Phi^k\big( z, H'(0), \1{H'(0)}, j_0^{k+1}H,
\1{j_0^{k+1}H} \big)
\end{equation}
holds for any local biholomorphism $H\colon (\bC^2,0)\to
(\bC^2,0)$ sending $M$ into $M'$, where \break
$\Phi^k(z,\L_1,\2\L_1,\L_2,\2\L_2)$ is a polynomial in
$(\L_2,\2\L_2)\in J_{0,0}^{k+1}(\C^2)\times
\1{J_{0,0}^{k+1}(\C^2)}$
with coefficients that are holomorphic
in $(z,\L_1,\2\L_1)\in
\C\times J_{0,0}^{1}(\C^2)\times \1{J_{0,0}^{1}(\C^2)}$
in a neighborhood of $\{0\}\times
\T(\bC^2)\times \1{\T(\bC^2)}$. Moreover, the $\bC^2$-valued
functions $\Phi^k$ depend only on $(M,p)$ and $(M',p')$.
\end{Thm}

Observe, that in normal coordinates, the line $\{(z,0) : z\in
\C\}$ is the zero Segre variety and hence we see the conclusion of
Theorem~\ref{thm-jetsegre} as a parametrization of jets along the
zero Segre variety. We first show how Theorem~\ref{thm-jetsegre}
implies Theorem~\ref{dynamics}.

\bpf[Proof of Theorem $\ref{dynamics}$] Let $H$ be a local
biholomorphism tangent to the identity at $0$, as in Theorem
\ref{dynamics}, and $(M,0)$ the germ of a Levi-nonflat
real-analytic hypersurface which is preserved by $H$. By
assumption, we have $j^1_0 H= j^1_0 \id$. Thus, by
Theorem~\ref{thm-jetsegre}, with $k=0$, applied to the local
biholomorphisms $H$ and $\id$, both sending $M$ into itself, we
conclude that $H(z,0)\equiv z$. Hence each point of the complex
hypersurface $\{z=0\}\subset\C^2$ is a fixed point of $H$. This
completes the proof of Theorem \ref{dynamics}. \epf

Before entering the proof of Theorem \ref{thm-jetsegre}, we shall
introduce some notation.
The equation (\ref{eq-normal1}) can be written in complex form
\begin{equation}\Label{eq-normal2}
M: \quad w=Q(z,\bar z,\bar w),
\end{equation}
where $Q(z,\chi,\tau)$ is a holomorphic function
satisfying
\begin{equation}\Label{eq-normal3}
Q(z,0,\tau) \equiv Q(0,\chi,\tau) \equiv \tau.
\end{equation}
Equation (\ref{eq-normal2}) defines a real hypersurface if and
only if (see \cite{BER})
\begin{equation}\Label{eq-normal4}
Q(z,\chi,\bar Q(\chi,z, w))\equiv w,
\end{equation} where we use the
notation $\bar h(\zeta):=\overline{ h(\bar \zeta)}$.

\renewcommand{\t}{\tau}

We shall study the derivatives of $Q(z,\chi,\tau)$ with respect
to $(z,\t)$ denoted as follows:
\begin{equation}\Label{eq-exp1}
q_{\alpha\mu}(\chi):=Q_{z^\a \t^\mu}(0,\chi,0), \quad
q_{\alpha\mu}(0)=0, \quad
\a\ge 1.
\end{equation}
Note that we do not use the function $q_{\a\mu}(\chi)$ for $\a=0$,
since the corresponding derivatives can be computed directly by
(\ref{eq-normal3}):
\begin{equation}\Label{Q-der-id}
Q_{\t^\mu}(0,\chi,0)\equiv
\begin{cases}
1 & \text{ for } \mu=1 \\
0 & \text{ for } \mu>1.
\end{cases}
\end{equation}
We now use the functions $q_{\a\mu}(\tau)$ to define biholomorphic
invariants of $(M,p)$ as follows. Let $m_0$ be the positive
integer (or $\infty$) given by
\begin{equation}\Label{eq-m0}
m_0:=\min \{m\in \bZ_+ : q_{\alpha\mu}(\chi)\not \equiv 0,\
\alpha+\mu=m\},
\end{equation}
where we set $m_0=\infty$ if $q_{\alpha\mu}(\chi)\equiv 0$ for all
$(\alpha,\mu)$. Thus, $m_0=\infty$ is equivalent to $M$ being
Levi-flat. In what follows, we shall assume that $M$ is not
Levi-flat, i.e. $m_0<\infty$. We also define
\begin{equation} \Label{eq-alpha0}
\mu_0:=\min\{\mu\in \bZ_+ : q_{\alpha\mu}(\chi)\not\equiv 0,\
\alpha+\mu=m_0\},
\end{equation}
and set $\alpha_0:=m_0-\mu_0$.

Let $M'\subset \bC^2$ be another real-analytic hypersurface, and
let $(z',w')\in\bC^2$ be normal coordinates for $M'$ at some point
$p'$. In what follows, we shall use a $'$ to denote an object
associated to $M'$ corresponding to one defined previously for
$M$. Let furthermore $H=(F,G)$ be a local mapping $(\bC^2,p)\to
(\bC^2,p')$. Then $H$ sends (a neighborhood of $p$ in) $M$ into
$M'$ if and only if it satisfies the identity
\begin{equation}\Label{eq-map1}
G(z,Q(z,\chi,\tau))\equiv
Q'\big(F(z,Q(z,\chi,\tau)),\bar F(\chi,\tau), \bar G(\chi,\tau)\big),
\end{equation}
for $(z,\chi,\tau)\in\bC^3$. In particular, setting $\chi=\t=0$,
we deduce that
\begin{equation}\Label{g-id}
G(z,0)\equiv 0.
\end{equation}
It follows from (\ref{g-id}) that, in normal coordinates, the
$2\times 2$ matrix $H'(0)$ is triangular and therefore $H$ is a
local biholomorphism if and only if
\begin{equation}\Label{eq-locbi}
F_z(0)\, G_w(0)\neq 0.
\end{equation}

We have the following property.

\begin{Pro}\Label{prop-m0} Let $M\subset  \bC^2$ be a real-analytic
hypersurface. Then, the integers $m_0$, $\alpha_0$, and $\mu_0$
defined above are biholomorphic invariants.
\end{Pro}

\begin{proof}
We have to show that, if $(M,p)$ and $(M',p')$ are locally
biholomorphic, then $\mu_0=\mu'_0$ and $\alpha_0=\alpha'_0$. We
introduce the following ordering of the pairs $(\alpha,\mu)\in
\bZ_+^2$. We write $(\alpha,\mu)\prec(\beta,\nu)$ if either
$\alpha+\mu<\beta+\nu$, or if $\alpha+\mu=\beta+\nu$ and $\mu<\nu$
(or, equivalently, $\alpha>\beta$). We prove the statement by
contradiction. Suppose $(\alpha_0,\mu_0)\prec(\alpha'_0,\mu'_0)$.
We differentiate (\ref{eq-map1}) by the chain rule $\a_0$ times
in $z$ and $\mu_0$ times in $\t$, evaluate the result at
$(z,\t)=0$ and use the identities (\ref{eq-normal3}),
(\ref{Q-der-id}) and (\ref{g-id}). On the right-hand side we
obtain a sum of terms each of which has a factor of $Q'_{z^\a
\chi^\b \t^\mu}(0,\bar F(\chi,0),0)$ with $(\a,\mu)\preceq
(\a_0,\mu_0)$. By the assumption
$(\alpha_0,\mu_0)\prec(\alpha'_0,\mu'_0)$, all these derivatives
are zero. Similarly, using the fact that $Q_{z^\a
\t^\mu}(0,\chi,0) \equiv 0$ for
$(\alpha,\mu)\prec(\alpha_0,\mu_0)$ on the left-hand side, we
conclude that
\begin{equation}\Label{eq-mu0-contra}
G_{z^{\a_0}w^{\mu_0}}(0) +
G_w(0)\,q_{\alpha_0\mu_0}(\chi) \equiv 0.
\end{equation}
Since $G_w(0)\neq 0$ and $q_{\alpha_0\mu_0}(\chi)\ne\const$,
we reach the desired contradiction.
We must then have $(\alpha_0,\mu_0)\succeq (\alpha'_0,\mu'_0)$.
The opposite inequality follows by reversing the roles of $M$ and
$M'$ by considering the inverse mapping $H^{-1}$. Hence,
$(\alpha_0,\mu_0)=(\alpha'_0,\mu'_0)$ as claimed.
\end{proof}

We are now ready to prove Theorem \ref{thm-jetsegre}.

\begin{proof}[Proof of Theorem $\ref{thm-jetsegre}$]
As in the proof of Proposition \ref{prop-m0},
we differentiate (\ref{eq-map1}) by the chain rule
$\a_0$ times in $z$ and $\mu_0$ times in $\t$,
evaluate the result at $(z,\t)=0$ and use the identities
(\ref{eq-normal3}),
(\ref{Q-der-id}) and (\ref{g-id}):
\begin{equation}\Label{eq-mu0-0}
G_{z^{\a_0}w^{\mu_0}}(0)+
G_w(0)\,q_{\alpha_0\mu_0}(\chi)\equiv
q'_{\alpha_0\mu_0}(\bar
F(\chi,0))\big(F_z(0)+F_w(0)
q_{10}(\chi)\big)^{\alpha_0}
\1{G_w}(\chi,0)^{\mu_0}.
\end{equation}
By putting $\chi=0$ and using (\ref{eq-exp1}) we obtain
$G_{z^{\a_0}w^{\mu_0}}(0)=0$. Furthermore, the derivative
$\1{G_w}(\chi,0)$ on the right-hand side of (\ref{eq-mu0-0}) can
be computed by differentiating (\ref{eq-map1}) in $\t$ at
$(z,\t)=0$:
\begin{equation}\Label{eq-Gw}
\1{G_w}(\chi,0)\equiv G_w(0) - q'_{10}(\bar F(\chi,0)) F_w(0).
\end{equation}
Note that for $\chi=0$, we obtain $\1{G_w}(0)= G_w(0)$, i.e.
$G_w(0)$ is real. After these observations, (\ref{eq-mu0-0}) can
be rewritten as
\begin{equation}\Label{eq-mu0}
G_w(0)\,q_{\alpha_0\mu_0}(\chi)\equiv q'_{\alpha_0\mu_0}(\bar
F(\chi,0)) \big(F_z(0)+ q_{10}(\chi)F_w(0)\big)^{\alpha_0}
\big({G_w}(0)- q'_{10}(\bar F(\chi,0))F_w(0)\big)^{\mu_0}.
\end{equation}

We point out that $q_{10}(\chi)\not \equiv 0$ (and similarly
$q'_{10}(\chi')\not\equiv 0$) if and only if $m_0=1$ (in this
case, $(\alpha_0,\mu_0)=(1,0)$). It can be shown that $m_0=1$ if
and only if $M$ is {\it finitely nondegenerate} at $p$ (see e.g.\
\cite{BER}). We will not use this fact in the present paper.

It follows, that for $m_0\geq 2$, (\ref{eq-mu0}) reduces to the form
\begin{equation}\Label{eq-mu0-mgeq2}
q_{\alpha_0\mu_0}(\chi)\equiv
q'_{\alpha_0\mu_0}(\bar F(\chi,0))F_z(0)^{\alpha_0}G_w(0)^{\mu_0-1}.
\end{equation}
If $m_0=1$, in which case $(\alpha_0,\mu_0)=(1,0)$, we instead have
\begin{equation}\Label{eq-mu0-m=1}
{G_w}(0)\,q_{10}(\chi)\equiv
q'_{10}(\bar F(\chi,0))
\big(F_z(0)+F_w(0) q_{10}(\chi)\big).
\end{equation}
Define $l\ge 1$ to be the order of
vanishing at $0$ of the function $q_{\alpha_0\mu_0}(\chi)\not\equiv 0$.
(It is
not difficult to verify from (\ref{eq-normal4}) that $l\geq \alpha_0$.
We will not use this property.)
Together with (\ref{eq-locbi}),
it follows from (\ref{eq-mu0-mgeq2}) and (\ref{eq-mu0-m=1}),
respectively, that $l=l'$. Let us write
\begin{equation}\Label{eq-factor1}
q_{\alpha_0\mu_0}(\chi)\equiv (q(\chi))^{l}, \quad
q'_{\alpha_0\mu_0}(\chi')\equiv (q'(\chi'))^{l},
\end{equation}
where $q$ and $q'$ are local
biholomorphisms of $\C$ at $0$.
Taking $l$th roots from both sides of (\ref{eq-mu0-mgeq2}) and
(\ref{eq-mu0-m=1}) respectively we obtain
\begin{equation}\Label{eq-Fdet}
q'(\bar F(\chi,0))\equiv E(\chi)\,q(\chi),
\end{equation}
where
\begin{equation*}
E(\chi):=
\begin{cases}
\Big(G_w(0)\big(F_z(0) + F_w(0)q_{10}(\chi)\big)^{-1}\Big)^{1/l}, &
m_0=1 \\
\big(F_z(0)^{-\alpha_0}G_w(0)^{1-\mu_0}\big)^{1/l}, & m_0\geq 2,
\end{cases}
\end{equation*}
and the branch of the $l$th root has to be chosen appropriately.
We want to write $E$ as a holomorphic function in $\chi$,
$F_z(0)$, $\1{F_z}(0)$ and $F_w(0)$. For this, observe that there
exists a constant $c\ne 0$ such that
\begin{equation}\Label{ansatz}
E(\chi)\equiv
\begin{cases}
c\big(1 + F_w(0)F_z(0)^{-1}q_{10}(\chi)\big)^{-1/l}, &  m_0=1 \\
c, & m_0\geq 2,
\end{cases}
\end{equation}
where we used the principal branch of the $l$th root near $1$
(for $\chi$ small). We substitute the expressions (\ref{ansatz})
for $E(\chi)$ in (\ref{eq-Fdet}) and differentiate once in $\chi$
at $\chi=0$ to obtain
\begin{equation}\Label{eq-c}
q'_{\chi'} (0) \1{F_z}(0) \equiv cq_\chi (0)
\end{equation}
in both cases $m_0=1$ and $m_0\ge 2$.
Solving (\ref{eq-c}) for $c$ and using (\ref{eq-Fdet}) and
(\ref{ansatz}) we obtain
\begin{equation}\Label{Fdet}
\bar F(\chi,0) \equiv \chi \Psi(\chi ,H'(0),\1{H'(0)}), \quad
\end{equation}
where $\Psi$ is a holomorphic function in its arguments, defined
in a neighborhood of the subset
$\{0\}\times \T(\bC^2)\times \1{\T(\bC^2)}
\subset \C\times J^1_{0,0}(\C^2)\times \1{J^1_{0,0}(\C^2)}$
such that
\begin{equation}\Label{nonvanish}
\Psi(0,\L, \2\L)\ne 0
\end{equation}
for any $(\L,\2\L)\in \T(\bC^2)\times \1{\T(\bC^2)}$.
In particular, we
see from (\ref{Fdet}) that the mapping $H$ along the Segre variety
$\{w=0\}$ is completely determined by its first jet at
$0$. In fact $H(z,0)$ depends only on the derivatives $F_w(0)$,
$F_z(0)$, $\1{F_z}(0)$.

We shall now determine derivatives of $H$ with respect to $w$
along the Segre variety $\{w=0\}$ in terms of jets of $H$ at the
origin. We begin with the derivatives $F_{w^k}(z,0)$ (or
$\1{F_{w^k}}(\chi,0)$),  $k\ge 1$. To get an expression involving
them we use the same strategy as above, but now we differentiate
(\ref{eq-map1}) $\a_0$ times in $z$ and $\mu_0+k$ times in $\t$
and then set $(z,\t)=0$. We shall use the notation $(F_{z^i
w^j})_{i+j\le k+1}$ etc. to denote strings of partial derivatives
of {\em positive order}. We obtain
\begin{multline}\Label{long}
G_{z^{\a_0} w^{\mu_0+k}}(0) + \Psi_1(\chi,j_0^{k+1}G) \equiv \\
\d_{\z} q'_{\a_0\mu_0}(\bar F(\chi,0)) \big(F_z(0) + F_w(0)q_{10}
(\z)\big)^{\a_0} \1{F_{w^k}}(\chi,0) \1{G_w}(\chi,0)^{\mu_0} +\\
\Psi_2\Big(\chi,\bar F(\chi,0), j_0^{k+1}F,
\big(\1{F_{w^r}}(\chi,0)\big)_{r\le k-1},
\big(\1{G_{w^s}}(\chi,0)\big)_{s\le k+1} \Big),
\end{multline}
where the functions $\Psi_1(\chi,\L)$ and
$\Psi_2(\chi,\chi',\L_1,\L_2,\L_3)$ are polynomials in $\L$ and
$(\L_1,\L_2,\L_3)$ respectively with holomorphic coefficients in
$(\chi,\chi')$. (More precisely, the coefficients are polynomials
in $(q_{\alpha\mu}(\chi))$ and $(q'_{\alpha\mu}(\chi'))$ and in
their derivatives.) We claim, however, that no term involving
$\1{G_{w^{k+1}}}(\chi,0)$ occurs with a nontrivial coefficient in
$\Psi_2$ when $m_0=1$. Indeed, in this case $\alpha_0=1$ and
$\mu_0=0$. Thus, to obtain the expression (\ref{long}) we
differentiate once in $z$ and $k$ times in $\tau$. Consequently,
no term of the form $\1{G_{w^{k+1}}}(\chi,0)$ can appear, as
claimed.

We observe that from the identity (\ref{long}) with $\chi=0$ we
obtain a polynomial expression for the derivative $G_{z^{\a_0}
w^{\mu_0+k}}(0)$ in terms of $j_0^{k+1}H$ and $\1{j_0^{k+1}H}$.
(Recall that $H=(F,G)$.) After substituting in (\ref{long}) this
expression for $G_{z^{\a_0} w^{\mu_0+k}}(0)$, the right-hand side
of (\ref{eq-Gw}) for $\1{G_w}(\chi,0)$ and the right-hand side of
(\ref{Fdet}) for $\bar F(\chi,0)$,  we obtain
\begin{multline}\Label{long1}
\d_{\chi'} q'_{\a_0\mu_0} \big(\chi\Psi(\chi,H'(0),\1{H'(0)}\big)
\, \1{F_{w^k}}(\chi,0) \equiv \\ \Psi_3\Big(\chi, H'(0),
\1{H'(0)}, j_0^{k+1}H, \1{j_0^{k+1}H},
\big(\1{F_{w^r}}(\chi,0)\big)_{r\le k-1},
\big(\1{G_{w^s}}(\chi,0)\big)_{s\le k+1} \Big),
\end{multline}
where $\Psi_3(\chi,\L_1,\2\L_1,\L_2,\2\L_2,\L_3,\L_4)$ is a
polynomial in  $(\L_2,\2\L_2,\L_3,\L_4)$ with holomorphic
coefficients in $(\chi,\L_1,\2\L_1)$. We also observe that the
coefficient of $\1{F_{w^k}}(\chi,0)$ on the left-hand side does
not vanish identically by the choice of $(\a_0,\mu_0)$.

As before, we get an expression for the derivatives $\1{G_{w^s}}
(\chi,0)$ that occur on the right-hand side of (\ref{long1}) by
differentiating (\ref{eq-map1}) in $\t$ at $(z,\t)=0$, this time
$s\ge 2$ times:
\begin{multline}\Label{eq-Gwk}
G_{w^s}(0)\equiv q'_{10}(\bar F(\chi,0)) F_{w^s}(0) +\\ \d_{\chi'}
q'_{10}(\bar F(\chi,0)) F_{w}(0) \1{F_{w^{s-1}}}(\chi,0)+
\1{G_{w^s}}(\chi,0)+\\ \Psi_4\Big( \bar F(\chi,0),
\big(F_{w^i}(0)\big)_{i\le s-1},
\big(\1{F_{w^r}}(\chi,0)\big)_{r\le s-2},
\big(\1{G_{w^j}}(\chi,0)\big)_{j\leq s-1} \Big),
\end{multline}
where $\Psi_4(\chi',\L_1,\L_2,\L_3)$ is a polynomial in
$(\L_1,\L_2,\L_3)$ with holomorphic coefficients in  $\chi'$. We
solve (\ref{eq-Gwk}) for $\1{G_{w^s}}(\z,0)$ in terms of
$\1{G_{w^{t}}}(\z,0)$ with $t<s$. Then, by induction on $s$ and by
(\ref{eq-Gw}), we obtain for any $s\ge 1$ an identity of the form
\begin{multline}\Label{eq-Gwk1}
\1{G_{w^s}}(\chi,0) \equiv G_{w^s}(0) -q'_{10}(\bar F(\chi,0))
F_{w^s} (0)- \d_{\chi'} q'_{10}(\bar F(\chi,0)) F_{w}(0)
\1{F_{w^{s-1}}}(\chi,0)+ \\ \Psi_5\Big( \bar F(\chi,0),
\big(H_{w^i}(0)\big)_{i\le s-1},
\big(\1{F_{w^r}}(\chi,0)\big)_{r\le s-2} \Big),
\end{multline}
where $\Psi_5(\chi',\L_1,\L_2)$ is a polynomial in $(\L_1,\L_2)$
with holomorphic coefficients in $\chi'$. Note that the term
$\1{F_{w^{s-1}}}(\chi,0)$ only occurs in  (\ref{eq-Gwk1}) when
$m_0=1$.

We now substitute the right-hand side of (\ref{eq-Gwk}) for
$\1{G_{w^s}}(\chi,0)$ and the right-hand side of (\ref{Fdet}) for
$\bar F(\chi,0)$ in the identity (\ref{long1}) to obtain
\begin{equation}\Label{eq-Fwk}
\d_{\chi'} q'_{\a_0\mu_0} \big(\chi\Psi(\chi,H'(0),\1{H'(0)}\big)
\, \1{F_{w^k}}(\chi,0)\equiv \Psi^k \Big(\chi, H'(0), \1{H'(0)},
j_0^{k+1}H, \1{j_0^{k+1}H}, \big(\1{F_{w^r}}(\chi,0)\big)_{r\le
k-1} \Big),
\end{equation}
where $\Psi^k(\chi,\L_1,\2\L_1,\L_2,\2\L_2,\L_3)$ a polynomial in
$(\L_2,\2\L_2,\L_3)$ with holomorphic coefficients in \break
$(\chi,\L_1,\2\L_1)$. We claim that an identity of the form
\begin{equation}\Label{eq-jetsegre-f}
\bar F_{w^k}(\chi,0) \equiv \2\Psi^k\big(\chi, H'(0), \1{H'(0)},
j_0^{k+1}H , \1{j_0^{k+1}H}\big)
\end{equation}
holds, where $\2\Psi^k(\chi,\L_1,\2\L_1,\L_2,\2\L_2)$ is a
polynomial in $(\L_2,\2\L_2)$ with holomorphic coefficients in
$(\chi,\L_1,\2\L_1)$. We prove the claim by induction on $k$. For
$k=1$, there is no occurrence of $\L_3=\1{F_{w^r}}(\chi,0)$ on the
right-hand side of (\ref{eq-Fwk}). We now would like to divide
both sides of (\ref{eq-Fwk}) by the first factor
$\Gamma:=\d_{\chi'} q'_{\a_0\mu_0}
\big(\chi\Psi(\chi,H'(0),\1{H'(0)}\big)$ on the left-hand side. We
know from (\ref{nonvanish}) that this factor does not vanish
identically. However, it may happen that it vanishes for $\chi=0$.
In order to obtain a holomorphic function on the right-hand side
after this division, we have to make sure that the vanishing
order of the right-hand side with respect to $\chi$ at the origin
is not smaller that the vanishing order of $\Gamma$. Of course,
by (\ref{eq-Fwk}), this is true for those jet values of
$(\L_1,\2\L_1,\L_2,\2\L_2)$ that come from a local biholomorphism
$H$ sending $M$ into $M'$. On the other hand, we have no
information about the vanishing order of
$\Psi^1(\chi,\L_1,\2\L_1,\L_2,\2\L_2)$ for other values of
$(\L_1,\2\L_1,\L_2,\2\L_2)$. Hence we may not be able to divide.
The idea for solving this problem is to extract the ``higher order
part'' of $\Psi^1$ that is divisible by $\Gamma$.

Let $0\le\nu<\infty$ be the vanishing order of the function
$\d_{\chi'} q'_{\a_0\mu_0}(\chi')$ at $\chi'=0$.
Then the vanishing order of the left-hand side of (\ref{eq-Fwk})
is at least $\nu$. By truncating the power series expansion in
$\chi$ of the coefficients of the polynomial $\Psi^1$ on the
right-hand side of (\ref{eq-Fwk}), we can write it in a unique
way as a sum $\Psi^1\equiv \Psi_1^1+\Psi_2^1$ such that both
$\Psi_1^1$ and $\Psi_2^1$ are polynomials in $(\L_2,\2\L_2)$ with
holomorphic coefficients in $(\chi,\L_1,\2\L_1)$, each coefficient
of $\Psi_1^1$ is a polynomial in $\chi$ of order at most $\nu-1$
with holomorphic coefficients in $(\L_1,\2\L_1)$ and each
coefficient of $\Psi_2^1$ is of vanishing order at least $\nu$
with respect to $\chi$. We now remark that, whenever we set
$(\L_1,\2\L_1,\L_2,\2\L_2)=(H'(0), \1{H'(0)}, j_0^{k+1}H,
\1{j_0^{k+1}H})$ for a mapping $H$ satisfying (\ref{eq-Fwk}), the
first polynomial $\Psi^1_1$ must vanish identically in $\chi$.
Therefore, the identity (\ref{eq-Fwk}) will still hold after we
replace $\Psi^1$ by $\Psi^1_2$. Now it follows from
(\ref{nonvanish}) that $\Psi^1_2$ is divisible by $\Gamma$ and,
hence, we have proved the claim for $k=1$. The induction step for
$k>1$ is essentially a repetition of the above argument. Once
(\ref{eq-jetsegre-f}) is shown for $k<k_0$, we substitute it for
$\1{F_{w^r}}(\chi,0)$ in the right-hand side of (\ref{eq-Fwk}).
Then we obtain a polynomial without $\L_3$ and the above argument
can be used to obtain (\ref{eq-jetsegre-f}) for $k=k_0$. The
claim is proved.

A formula for $\bar G_{w^k}(\chi,0)$ similar to
(\ref{eq-jetsegre-f}) is obtained by substituting the result for
$F$ in (\ref{eq-Gwk1}) (and also using (\ref{Fdet})). The proof
of Theorem~\ref{thm-jetsegre} is complete.
\end{proof}

\section{Parametrization of biholomorphisms in the finite type
  case}\Label{sec-finite}
In this section, we will prove the
  following theorem, from which Theorem \ref{main-jetpar} is a direct
consequence and Corollary
  \ref{thm-finite2} follows by standard techniques (see \cite{Asian}
and \cite{MA}). We keep the setup and notation
  introduced in previous sections.
We also denote by $\T^2(\C^2)$
the subspace of $2$-jets in $J^2_{0,0}(\C^2)$
whose first derivative matrix is upper triangular.

\begin{Thm}\Label{thm-jetpar}
Let $M,M'\subset\bC^2$ be real-analytic hypersurfaces of finite
type, and let  $(z,w)\in\bC^2$ and $(z',w')\in \bC^2$ be normal
coordinates for $M$ and $M'$ at $p\in M$ and  $p'\in M'$,
respectively. Then an identity of the form
\begin{equation*}
H(z,w) \equiv \Theta\big( z, w, j_0^2 H, \1{j_0^2 H} \big),
\end{equation*}
holds for any local biholomorphism $H\colon (\bC^2,p)\to
(\bC^2,p')$ sending $M$ into $M'$, where $\Theta(z,w,\L,\2\L)$
is a holomorphic function in  a neighborhood of the subset
$\{(0,0)\}\times \T^2(\bC^2)\times
\1{\T^2(\bC^2)}$
in $\C^2\times J^2_{0,0}(\C^2)\times \1{J^2_{0,0}(\C^2)}$,
depending only on $(M,p)$ and $(M',p')$.
\end{Thm}

\begin{proof}
We use the expansion of   the function   $Q(z,\chi,\tau)$ (and
similarly for the
function
  $Q'(z',\chi',\tau')$, using a $'$ to denote corresponding objects
  associated to $M'$) as follows
\begin{equation}\Label{eq-exp2}
  Q(z,\chi,\tau)=\tau+\sum_{\beta\ge 1} r_{\beta}(\chi,\tau)z^\beta,
\quad r_\b(0,0)=0.
\end{equation}
Recall that $M$ is of
  finite type at $0$ if and only if $Q(z,\chi,0)\not\equiv 0$,
i.e. if $r_\b(\chi,0)\not\equiv 0$ for some $\b$. We define a
positive integer associated to $(M,p)$ by
\begin{equation}\Label{eq-beta0}
\beta_0:=\min\{\beta\colon  r_\b(\chi,0)\not\equiv 0\}
\end{equation}
and an integer $\b'_0$ associated to $(M',p')$ by the analogous
formula.
 It follows easily, by setting $\tau=0$
  in the equation  (\ref{eq-map1}), differentiating in $z$ and
  then setting
 $z=0$ as in the proof of Proposition~\ref{prop-m0}, that
  $\beta_0=\beta'_0$
(i.e. $\beta_0$ is a biholomorphic invariant). Indeed, by using
also the fact that $G(z,0)\equiv 0$, we obtain in this way the
identity
\begin{equation}\Label{order1}
r_{\b_0}(\chi,0) G_w(0) \equiv  r'_{\b_0}(\bar F(\chi,0),0)
\big(F_z(0)  + F_w(0) r_1(\chi,0) \big).
\end{equation}
By differentiating (\ref{eq-map1})
  with respect to $z$ and setting $\tau=\bar Q(\chi,z,0)$, we also
  obtain
\begin{multline}\Label{eq-map1_z1} Q_z\big(z,\chi,\bar
  Q(\chi,z,0)\big)G_w(z,0)\equiv
\\ Q'_z\big(F(z,0),\bar
  F(\chi,\bar Q(\chi,z,0)),\bar G(\chi,\bar Q(\chi,z,0))\big) \;
  \big(F_z(z,0)+
  Q_z(z,\chi,\bar Q(\chi,z,0))F_w(z,0)\big),
\end{multline}
where we
  have
  used (\ref{eq-normal4}) and the fact that $G(z,0)\equiv 0$. By
  using the conjugate of (\ref{eq-map1}) to substitute for $\bar
G(\chi,\bar
  Q(\chi,z,0))$, we obtain
\begin{multline}\Label{eq-map1_z2}
Q_z\big(z,\chi,\bar Q\big)G_w(z,0)\equiv
\\Q'_z\big(F(z,0),\bar
  F(\chi,\bar Q),\bar Q'(\bar F(\chi,\bar
  Q),F(z,0), 0)\big) \;
\big(F_z(z,0)+
  Q_z(z,\chi,\bar Q)F_w(z,0)\big),
\end{multline}
where we have
  used the notation $\bar Q:=\bar Q(\chi,z,0)$. Observe, by
  differentiating (\ref{eq-normal4}) with respect to $z$ and
  setting $w=0$, that
\begin{equation}\Label{eq-normal-diff}
  Q_z\big(z,\chi,\bar Q(\chi,z,0)\big)\equiv
-Q_\tau\big(z,\chi,\bar Q(\chi,z,0)\big)   \bar  Q_z(\chi,z,0).
\end{equation}
Since $Q_\tau(z,0,0)\equiv 1$ by
  (\ref{eq-normal3}), we conclude from (\ref{eq-exp2}) and
  (\ref{eq-normal-diff}) that
  \begin{equation}\Label{eq-exp3}
Q_z\big(z,\chi,\bar Q(\chi,z,0)\big)
=   \d_z\1{r_{\beta_0}}(z,0) \chi^{\beta_0}+O(\chi^{\beta_0+1}).
  \end{equation}
  Also, observe that $\d_z\1{r_{\beta_0}}(z,0)\not\equiv 0$ by the
choice of
$\b_0$.
It follows that  $\d_z \1{r_{\beta_0}}(z,0)\neq 0$ for all
  $z$ in any sufficiently small punctured disc centered at 0.
In such a punctured disc there exist $\beta_0$ locally defined
holomorphic functions $\psi(z,\chi)$,
  differing by multiplication by a $\b_0$th root of unity, such that
  \begin{equation}\Label {eq-r-root}
  Q_z\big(z,\chi,\bar Q(\chi,z,0)\big)
\equiv \psi(z,\chi)^{\beta_0}
\end{equation}
and
\begin{equation}\Label{eq-r-root1}
\quad \psi(z,0)\equiv 0, \quad \psi_\chi(z,0)\neq 0
\text{ for } z\ne 0 \text{ near } 0
\end{equation}
for each choice of $\psi$. Moreover, each local function $\psi$
extends as a multiple-valued holomorphic function (with at most
$\beta_0$ branches) to a neighborhood of $(0,0)$ in
$\C^2\setminus\{z= 0\}$. If we choose $z$ in a sufficiently small
punctured disc centered at $0$,
  then also $\d_{z'} \1{r'_{\beta_0}}(F(z,0),0)  \neq  0$
and the equation (\ref{eq-map1_z2}) can be written
\begin{equation}\Label{eq-map1_z3}
  \psi(z,\chi)^{\beta_0}G_w(z,0) \big(F_z(z,0)+
Q_z(z,\chi,\bar Q)  F_w(z,0)\big)^{-1} \equiv
 \psi'  \big(F(z,0),\bar F(\chi,\bar Q(\chi,z,0))\big)^{\beta_0}.
\end{equation}
$G_w(z_0,0)\neq 0$.
Taking $\b_0$th roots from both sides of (\ref{eq-map1_z3})
and substituting the conjugate of (\ref{Fdet}) for $F(z,0)$ we obtain
\begin{equation}\Label{psi-id}
\psi' \big(z\bar\Psi(z,\1{H'(0)},H'(0)),\bar F(\chi,\bar
Q(\chi,z,0)) \big) \equiv E(z,\chi) \psi(z,\chi)
\end{equation}
with
\begin{equation*}
E(z,\chi):=\big(G_w(z,0)(F_z(z,0)+
\psi(z,\chi)^{\beta_0}F_w(z,0))^{-1} \big)^{1/\beta_0},
\end{equation*}
and the branch of the $\b_0$th root has to be appropriately
chosen. Analogously to the proof of Theorem~\ref{thm-jetsegre},
we write
\begin{equation}\Label{E-exp}
E(z,\chi)\equiv c \bigg(\frac{G_w(z,0)}{G_w(0)}
\Big(\frac{F_z(z,0)}{F_z(0)}+ Q_z(z,\chi,\bar Q)
\frac{F_w(z,0)}{F_z(0)} \Big)^{-1} \bigg)^{1/\beta_0},
\end{equation}
where $c$ is a constant and where we used the principal branch of the
$\b_0$th
root near $1$.
We now use Theorem~\ref{thm-jetsegre} to rewrite (\ref{E-exp}) as
\begin{equation}\Label{E-subs}
E(z,\chi)\equiv c\,
\Phi\big(z,\chi,H'(0),\1{H'(0)},j_0^2H,\1{j_0^2H}\big)^{1/\b_0},
\end{equation}
where $\Phi(z,\chi,\L_1,\2\L_1,\L_2,\2\L_2)$ is a holomorphic
function as in Theorem~\ref{thm-jetsegre}. Here the problem arises
that the values of $\Phi(0,0,\L_1,\2\L_1,\L_2,\2\L_2)$ may differ
from $1$ and thus the values of the $\b_0$th root in
(\ref{E-subs}) are not uniquely determined. To solve this problem
we use a trick to replace the function
$\Phi(z,\chi,\L_1,\2\L_1,\L_2,\2\L_2)$ by
$$\Phi(z,\chi,\L_1,\2\L_1,\L_2,\2\L_2) -
\Phi(0,0,\L_1,\2\L_1,\L_2,\2\L_2) + 1.$$ The new function (denote
it by $\Phi$ instead of the old one) takes the right value $1$ for
$(z,\chi)=0$, so that the principal  branch of its $\b_0$th root
is defined,  and the identity (\ref{E-subs}) still holds for any
local biholomorphism $H=(F,G)$ of $\C^2$ sending $(M,p)$ into
$(M',p')$.

We substitute the expression (\ref{E-subs}) for $E(z,\chi)$ in
(\ref{psi-id}) and differentiate once in $\chi$ at $\chi=0$ to
obtain
\begin{equation}\Label{eq-c2}
\psi'_{\chi'} \big(z\bar\Psi(z,\1{H'(0)},H'(0)),0\big) \equiv c\,
\Phi\big(z,0,H'(0),\1{H'(0)},j_0^2H,\1{j_0^2H}\big)^{1/\b_0}
\psi_\chi (z,0).
\end{equation}
By (\ref{eq-r-root}), the coefficient on the right-hand side does
not vanish identically. Hence we can solve $c$ from (\ref{eq-c2})
as a {\em multiple-valued holomorphic function}
\begin{equation}\Label{c}
c=c\big(z,H'(0),\1{H'(0)},j_0^2H,\1{j_0^2H}\big)
\end{equation}
and thus ignore the fact that $c$ is constant. We know from
(\ref{E-subs}), however, that, whenever the arguments of $\Psi$
and $\Phi$ are jets of a local biholomorphism, the ``function''
$c$ is actually (locally) constant. The different values of $c$
are due to the choice of different branches of $\psi$ and $\psi'$.
Recall that the values of $\psi$ and $\psi'$ may only differ by
$\b_0$th roots of unity. If $\psi$ and $\psi'$ are multiplied by
the roots of unity $\epsilon$ and $\epsilon'$ respectively, then
$c$ is multiplied by $\epsilon'\epsilon^{-1}$.

We claim that, for every fixed jet
$$(\L_1^0,\2\L_1^0,\L_2^0,\2\L_2^0)\in \T(\C^2)\times\1{\T(\C^2)}
\times J_{0,0}^2(\C^2)\times\1{J_{0,0}^2(\C^2)},$$ the
(multiple-valued) function $c(z,\L_1,\2\L_1,\L_2,\2\L_2)$ is
uniformly bounded for $(z,\L_1,\2\L_1,\L_2,\2\L_2)$ in some
neighborhood of $(0,\L_1^0,\2\L_1^0,\L_2^0,\2\L_2^0)$. Indeed, it
follows from the construction that the derivative $\psi_\chi(z,0)$
equals the $\b_0$th root of $\d_z \1{r_{\b_0}}(z,0)$. Furthermore,
by equation (\ref{order1}), the functions $\d_z \1{r_{\b_0}}(z,0)$
and $\d_z \1{r'_{\b_0}}(F(z,0),0)$ have the same vanishing order
at $z=0$. Then, taking $\b_0$th powers of both sides in
(\ref{eq-c2}) we obtain (single-valued) holomorphic functions of
the same vanishing order at $z=0$ (recall that both functions
$\Psi$ and $\Phi$ do not vanish at $z=0$). Hence $c^{\b_0}$ is
bounded as a ratio of two holomorphic functions having the same
vanishing order. The claim is proved.

We now substitute (\ref{c}) for $c$ into (\ref{E-subs}) and then use
(\ref{psi-id})  to obtain
\begin{multline}
  \Label{eq-map1_z4}
  \psi'\big(z\bar\Psi(z,\1{H'(0)},H'(0)),\bar F(\chi,\bar
  Q(\chi,z,0))\big)\equiv \\
c\big(z,H'(0),\1{H'(0)},j_0^2H,\1{j_0^2H}\big) \;
\Phi\big(z,\chi,H'(0),\1{H'(0)},j_0^2H,\1{j_0^2H}\big)^{1/\b_0}
\psi(z,\chi)
\end{multline}
 Since $\psi'_\chi(F(z,0),0)\neq 0$ for $z\ne 0$ near $0$, we may
 apply the implicit
  function theorem and solve for $\bar F(\chi,\bar Q(\chi,z,0))$ in
  (\ref{eq-map1_z4})
to conclude that
  \begin{equation}\Label{eq-barFzQ}
\bar F(\chi,\bar Q(\chi,z,0)) \equiv
\Phi_1(z,\chi,j_0^2H,\1{j_0^2H}),
  \end{equation}
where $\Phi_1(z,\chi,\L,\2\L)$ is an (a priori multiple-valued)
holomorphic function defined in a domain
\begin{equation}\Label{domain}
D\subset \C\times\C\times J_{0,0}^2(\C^2)\times\1{J_{0,0}^2(\C^2)}
\end{equation}
that contains all points $(z_0,0,\L_0,\2\L_0)$  with $z_0\ne 0$
that are sufficiently close to the set
$\{0\}\times\{0\}\times\T^2(\C^2)\times\1{\T^2(\C^2)}$. A value
of $\Phi_1$ depends on the values of $\psi$, $\psi'$ and $c$. We
have seen that, if $\psi$ and $\psi'$ are multiplied by $\epsilon$
and $\epsilon'$ respectively, then $c$ is multiplied by
$\epsilon'\epsilon^{-1}$. We now observe, that in fact the
identity (\ref{eq-map1_z4}) is  invariant under this change. Hence
we obtain {\em exactly the same value} for  $\Phi_1$ for {\em all
possible values} of $\psi$ and $\psi'$. We conclude that $\Phi_1$
is single-valued.

A similar expression for $\bar G$ is obtained
by substituting (\ref{eq-barFzQ}) and (\ref{Fdet}) into
(\ref{eq-map1}).
In summary, we obtain
\begin{equation}\Label{eq-HzQ}
\bar H(\chi,\bar Q(\chi,z,0)) \equiv \Xi\big(z,\chi,
j^2H(0),\1{j^2H(0)}\big),
\end{equation}
where $\Xi(z,\chi,\L,\2\L)$ is a ($\bC^2$-valued) holomorphic
function in a domain $D$ as in (\ref{domain}). Here we conjugated
and switched the variables.

To complete the proof of Theorem \ref{thm-jetpar}, we shall
proceed  as in \cite{Asian}. We consider the equation
\begin{equation}\Label{eq-w=Q1}
  \t= \bar Q(\chi,z,0),
\end{equation}
and try to solve it for $z$ as a function of $(\chi,\t)$ in a
neighborhood of a point $z_0\ne 0$ close to $0$. (We cannot apply
the method of \cite{Asian} at $0$, since the function $\Xi$ in
(\ref{eq-HzQ}) may not be defined in a neighborhood of
$\{0\}\times\{0\}\times \T^2(\C^2)\times\1{\T^2(\C^2)}$.) We
expand $\bar Q(\chi,z,0)$ in powers of $z-z_0$ near $z= z_0$ as
follows
\begin{equation}\Label{eq-exp4}
\bar Q(\chi,z,0)\equiv p_0(\chi;z_0)+ \sum_{\gamma\ge 1} p_\gamma
(\chi;z_0) (z-z_0)^\gamma
\end{equation}
with $p_\gamma(0;z_0) \equiv 0$ for all $\gamma$, since $\bar
Q(0,z,0)\equiv 0$. Let $\gamma_0$ be the smallest integer $\gamma
\geq 1$ such that $p_\gamma(\chi;z_0)\not\equiv 0$ with $z_0$
fixed. The existence of $\gamma_0<\infty$ is guaranteed by the
finite type  condition. Moreover, if $z_0\ne 0$ is sufficiently
small, $\gamma_0$ does not depend on $z_0$. After dividing
(\ref{eq-w=Q1}) by $p_{\gamma_0}(\chi;z_0)^{\gamma_0+1}$ and using
(\ref{eq-exp4}) we obtain
\begin{equation}\Label{eq-w=Q2}
  \frac{\t-p_0(\chi;z_0)}{p_{\gamma_0}(\chi;z_0)^{\gamma_0+1}} =
  \Big(\frac{z-z_0}{p_{\gamma_0}(\chi;z_0)}\Big)^{\gamma_0}+
\sum_{\gamma\ge\gamma_0+1} C_\gamma (\chi;z_0)
  \Big(\frac{z-z_0}{p_{\gamma_0}(\chi;z_0)}\Big)^{\gamma},
\end{equation}
where
$C_\gamma(\chi;z_0):=p_\gamma(\chi;z_0)p_{\gamma_0}(\chi;z_0)
^{\gamma-\gamma_0-1}$. Then, by the implicit function theorem, the
equation
\begin{equation*}
\eta=t^{\gamma_0 } + \sum_{\gamma\ge\gamma_0+1} C_\gamma
(\chi;z_0) t^{\gamma}
\end{equation*}
has $\gamma_0$ solutions of the form
$t=g(\chi,\eta^{1/\gamma_0};z_0)$, where $g(\chi,\zeta;z_0)$ is a
holomorphic function in a neighborhood of
$\{0\}\times\{0\}\times(\Delta_\delta\setminus\{0\})$ with
$g(0,0;z_0)\equiv 0$, where $\Delta_\delta:=\{z_0\in\C :
|z_0|<\delta\}$ is a sufficiently small disc. Hence the  equation
(\ref{eq-w=Q1}) can be solved for $z$ in the form
\begin{equation}\Label{eq-chi}
z= z_0+p_{\gamma_0}(\chi;z_0)\, g\bigg(\chi,
\Big(\frac{\t-p_0(\chi;z_0)} {p_{\gamma_0}(\chi;z_0)^{\gamma_0+1}}
\Big) ^{1/\gamma_0};z_0\bigg)
\end{equation}
for $\chi\ne 0$ and
$(\t-p_0(\z;z_0))/p_{\gamma_0}(\z;z_0)^{\gamma_0+1}$ both
sufficiently small.
We now substitute (\ref{eq-chi}) for $z$
in the identity  (\ref{eq-HzQ}) to obtain
\begin{equation}\Label{eq-HzQ3}
  \bar H(\chi,\t) \equiv \2\Xi\bigg(
\Big(\frac{\t-p_0(\chi;z_0)}{p_{\gamma_0}(\chi;z_0)^{\gamma_0+1}}
\Big)^{1/\gamma_0},\chi,
  j_0^2 H,  \1{j_0^2 H}; z_0
\bigg),
\end{equation}
where $\2\Xi(\zeta,\chi,\L,\2\L; z_0)$ is a holomorphic function
defined for all $(\zeta,\chi,\L,\2\L; z_0)$ with
$(\chi,\zeta;z_0)$ in a neighborhood of
$\{0\}\times\{0\}\times(\Delta_\delta\setminus\{0\})$ and
$$\big(z_0+p_{\gamma_0}(\chi;z_0)^{\gamma_0+1}
g(\chi,\zeta;z_0),\chi,\L,\2\L\big) \in D,$$ where $D$ is the
domain of definition of the function $\Xi$ in (\ref{eq-HzQ}). It
follows from the above description of $D$ that
$(0,0,\L_0,\2\L_0;z_0)$ is in the domain of definition of $\2\Xi$
whenever $z_0\ne 0$ and $(z_0,0,\L_0,\2\L_0)$ is in a
sufficiently small neighborhood of
$\{0\}\times\{0\}\times\T^2(\C^2)\times \1{\T^2(\C^2)}$.

 Let us expand the
function $\2\Xi(\zeta,\chi,\L,\2\L;z_0)$ in $\zeta$
\begin{equation}\Label{eq-decomp1}
\2\Xi(\zeta,\chi,\L,\2\L;z_0) \equiv \sum_{k\ge 0}
A_k(\chi,\L,\2\L;z_0) \zeta^k,
\end{equation} and
decompose it as $\2\Xi=\2\Xi_1+\2\Xi_2$, where
\begin{equation}\Label{eq-decomp2}
\begin{aligned}
    \2\Xi_1(\zeta,\chi,\L,\2\L;z_0) :=
&\,\sum_{j\ge 0} A_{j\gamma_0}(\chi,\L,\2\L;z_0) \zeta^{j\gamma_0} \\
    \2\Xi_2(\zeta,\chi,\L,\2\L;z_0) :=
&\,\sum_{{k\not \in \gamma_0\mathbb Z_+ }}
    A_k(\chi,\L,\2\L;z_0) \zeta^k.
\end{aligned}
\end{equation}
Since $\bar H(\chi,\t)$ is
    holomorphic in a neighborhood of $0$ in $\bC^2$, the function of
    $(\chi,\t;z_0)$ on
    the right hand side of (\ref{eq-HzQ3}) is independent of the
    value
of the $\gamma_0$th root, independent of $z_0\ne 0$ and extends
holomorphically
to a neighborhood of $0$ in $\C^3$.
Let us denote by
    $\L_0\in \T^2(\bC^2)$ the value of $j_0^2 H$. Since the
    function
    $\2\Xi_1(\zeta,\chi,\L_0,\bar \L_0;z_0)$, in which we substitute
\begin{equation}\Label{eq-s}
\zeta=\Big(\frac{\w-p_0(\chi;z_0)}{p_{\gamma_0}(\chi;z_0)^{\gamma_0+1}}
\Big)^{1/\gamma_0},
\end{equation}
is single valued on  $|\w|=\varepsilon$, for $\varepsilon>0$
sufficiently small  (depending on $\chi$ with
$p_{\gamma_0}(\chi;z_0) \neq 0$), we conclude that the
  function $(\zeta,\chi;z_0)\mapsto \2\Xi_2(\zeta,\chi,\L_0,\bar
  \L_0;z_0)$ is identically $0$.
Hence, we must have
\begin{equation}\Label{eq-HzQ4}
  \bar H(\chi,\t) \equiv \Xi_1
\Big(\frac{\w-p_0(\chi;z_0)}{p_{\gamma_0}(\chi;z_0)^{\gamma_0+1}},\chi,
  j_0^2 H,  \1{j_0^2 H};z_0\Big),
\end{equation}
where $\Xi_1(\eta,\chi,\L,\2\L;z_0) := \sum_{j\ge 0}
A_{j\gamma_0}(\chi,\L,\2\L;z_0) \eta^j$.

Next, we decompose
    each $A_{j\gamma_0}(\chi,\L,\2\L;z_0)$ uniquely as follows
\begin{equation}\Label{eq-Ajdecomp}
A_{j\gamma_0}(\chi,\L,\2\L;z_0)\equiv B_j(\chi,\L,\2\L;z_0)
p_{\gamma_0}(\chi;z_0)^{j(\gamma_0+1)}
    +\sum_{0\le l\le Kj-1}     R_{jl}(\L,\2\L;z_0)\chi^l,
\end{equation}
where $K$ denotes the order of
    vanishing of $p_{\gamma_0}(z;z_0)^{\gamma_0+1}$ at 0. It is not
    difficult to see that we have
\begin{equation}\Label{eq-Bj}
\sup_{|\chi|\leq \delta}|\!|B_j(\chi,\L,\2\L;z_0)|\!|\leq C^j
    \sup_{|\chi|\leq \delta} |\!|A_{j\gamma_0}(\chi,\L,\2\L;z_0)|\!|,
\end{equation} for
    some small $\delta>0$ and constant $C$, where $|\!| v|\!|$
    denotes
    the maximum of $|v_1|$ and $|v_2|$ for $v\in \bC^2$. Hence, the
    power series
\begin{equation}\Label{eq-Theta}
    \Gamma(\k,\chi,\L,\2\L;z_0):=\sum_{j\ge 0} B_j(\chi,\L,\2\L;z_0)
    \k^j
\end{equation}
defines a holomorphic function whose domain of definition
contains any point
$(0,0,\L_0,\2\L_0;z_0)$ with $z_0\ne 0$ and $(z_0,0,\L_0,\2\L_0)$
in a sufficiently small neighborhood of $\{0\}\times\{0\}\times
\T(\C^2)\times \1{\T(\C^2)}$.

We now wish to show that $\Xi_1$ in (\ref{eq-HzQ4}) can be
replaced by $\Gamma$ with $\k=\t-p_0(\chi;z_0)$. For this, we
decompose the function $\Xi_1$ uniquely as
    $\Xi_1=\Xi_2+\Xi_3$, where
\begin{equation}\Label{eq-decomp3}\begin{aligned}
    \Xi_2(\eta,\chi,\L,\2\L;z_0) := &\,\sum_{j\ge 0}
    B_{j}(\chi,\L,\2\L;z_0)p_{\gamma_0}(\chi;z_0)^{j(\gamma_0+1)}
    \eta^{j} \\
\Xi_3(\eta,\chi,\L,\2\L;z_0) := &\,\sum_{j\ge 0}
R_{j}(\chi,\L,\2\L;z_0) \eta^{j},
\end{aligned}
\end{equation} where
    $R_j(\chi,\L,\2\L;z_0):=\sum_{0\le l\le Kj-1} R_{jl}(\L,\2\L;z_0)
    \chi^l$ is the remainder polynomial in the division
    (\ref{eq-Ajdecomp}).
Now, observe that
    $\Xi_2(\eta,\chi,\L,\2\L;z_0)$,
with $\eta=\k/{p_{\gamma_0}(\chi;z_0)^{\gamma_0+1}}$, coincides
with the function $\Gamma(\k,\chi,\L,\2\L;z_0)$. Since the
right-hand side of (\ref{eq-HzQ4}) is holomorphic in $(\chi,\t)$
near $0$,
 it is  not difficult to see that $(z,\zeta)\mapsto \Xi_3(\eta,\chi,
 \L_0,\bar \L_0;z_0)$
must be identically $0$, and that
\begin{equation}\Label{eq-HzQ5}
  \bar H(\chi,\t) \equiv \Gamma\big(\t-p_0(\chi;z_0),\chi, j_0^2 H,
  \1{j_0^2 H};z_0 )\big).
\end{equation}
It remains to remark that the right-hand side of (\ref{eq-HzQ5})
is  holomorphic in $(\chi,\t;z_0)$ near $(0,0,\2{z_0})$ with any
$\2{z_0}\ne 0$ sufficiently small and is independent of $z_0$. The
proof of Theorem \ref{thm-jetpar} is complete.
\end{proof}

\bpf[Proof of Corollary~\ref{thm-finite2}]
The proof can be obtained by repeating the arguments from \cite{Asian}.
\epf

\renewcommand{\t}{\theta}
\section{Finite jet determination for solutions of singular
ODEs}\Label{sec-ode}

The proof of Theorem \ref{main}
in the infinite type case is based on Theorem~\ref{thm-jetsegre},
on the first author's results in \cite{Enonm}
(see Theorem~\ref{E-cor} below)
and on the following property of solutions of singular
ordinary differential equations which we prove in this section,
and which may be of independent interest.

\begin{Thm}\Label{ode}
Consider a singular differential equation for an $\bR^n$-valued
function $y(x,\theta)$, where $x\in\bR$, $\theta\in \bR^m$, of the form
\begin{equation}\Label{eq-m}
x^{\gamma+1} \partial_x y(x,\t)=
\frac{p(x,y(x,\t),\t)}{q(x,y(x,\t),\t)},
\end{equation}
where $\gamma\ge 0$ is an integer,
$p(x,y,\t)$ and $q(x,y,\t)$ are real-analytic functions
(valued in $\R^n$ and $\R$, respectively)
defined in a neighborhood of $0$ in $\R\times\R^n\times\R^m$
with $q(0,0,\t)\not\equiv 0$.
Let $\hat y(x,\t)$ be a real-analytic solution of $(\ref{eq-m})$
near $0$ with $\hat y(0,\t) \equiv 0$.
Then there exists an integer $k\geq 0$
such that, if $y(x,\t)$ is another solution near the origin
with $\d_x^l y(0,\t) \equiv \d_x^l \hat y(0,\t)$ for $0\le l\le k$,
then $y(x,\t)\equiv \hat y(x,\t)$.
\end{Thm}

\begin{proof}
We write $f(x,y,\t)$ for the right-hand side of (\ref{eq-m}).
For $\gamma=0$, the proof is rather simple.
Expand both sides of the equation (\ref{eq-m}) in
  powers of $x$ and identify the coefficients of $x^k$. It is not
 difficult to see that one may solve the resulting equation for
 the coefficients $a_k=a_k(\t)$ of $y(x,\t)$ in terms of $a_l=a_l(\t)$,
 with  $l\leq k-1$, unless $k$ is an
 eigenvalue of $f_y(0,0,\t)$.
The last possibility does not happen if $k$ is
sufficiently large and $\t$ is outside a countable union $\Gamma$
of proper real-analytic subvarieties.
For every $\t\notin\Gamma$ and $y,\hat y$ satisfying the hypotheses in
the theorem,
we conclude that $y(x,\t)\equiv \hat y(x,\t)$.
The required statement follows by continuity and the fact
that the complement of $\Gamma$ is dense in a neighborhood of $0$ in
$\R^m$.
The details are left to the reader.

For the rest of the proof we assume $\gamma\ge 1$.
We shall write $y^{(s)}\in \R^n$ for the $s$th derivative of $y(x,\t)$
in $x$
evaluated at $(0,\t)$
and $y^{(s+1,\ldots,s+l)}\in \R^{ln}$ for the column of the $l$
derivatives
$y^{(s+1)},\ldots,y^{(s+l)}$.

We shall differentiate (\ref{eq-m}) in $x$ at $(0,\t)$
and apply the chain rule. Since $y(0,\t) \equiv 0$, the derivatives of
the right-hand side $f(x,y,\t)$ will be always evaluated at $(0,0,\t)$.
Hence they will be ratios of real-analytic functions where the
denominator is some power
of $q(0,0,\t)$. We consider the ring
of all ratios of this kind and all polynomials (and rational functions)
below
will be understood over this ring (i.e.\ a polynomial below will be a
polynomial with coefficients in this ring).

By taking the $s$th derivative ($s\ge
\gamma+1$) of the identity (\ref{eq-m}) in $x$, evaluating at $(0,\t)$,
and
using the chain rule, we obtain
\begin{equation}\Label{der-s}
c_s y^{(s-\gamma)} = P_l(y^{(1,\ldots,l-1)},\t) y^{(s-l+1,\ldots,s)}
+ R_{l,s}(y^{(1,\ldots,s-l)},\t)
\end{equation}
for any $1\le l \le s/2$, where $c_s=\binom{s}{ \gamma+1}$, $P_l$ is a
$\R^{n\times nl}$-valued matrix polynomial in $y^{(1,\ldots,l-1)}$
depending only on $l$ and $R_{l,s}$ is a $\R^n$-valued polynomial
in $y^{(1,\ldots,s-l)}$ depending on both $l$ and $s$.
In fact, $P_l$ can be written as $P_l=(P_l^1,\ldots,P_l^l)$
with each $P_l^i$ being an $n\times n$ matrix
given by the formula
$$P_l^i(y^{(1,\ldots,l-1)},\t) =  (d/dx)^{l-i}(f_y(x,y(x),\t))|_{x=0},
\quad i=1,\ldots,l.$$

We further fix integers $t$ and $r$
satisfying $t+1 \le (r+1/\gamma)/2$,
collect the identities (\ref{der-s}) in blocks for
$l=(t+1)\gamma$ and $s=r\gamma+1,\ldots,(r+1)\gamma$
and write them in the form
\begin{equation}\Label{collect}
C_r y^{((r-1)\gamma+1,\ldots,r\gamma)} =
\sum_{0\le j\le t} Q^j
(y^{(1,\ldots,(j+1)\gamma-1)},\t) y^{((r-j)\gamma+1,\ldots,
(r-j+1)\gamma)} + S_{r,t}(y^{(1,\ldots,(r-t)\gamma)},\t),
\end{equation}
where $C_r$ is the diagonal $\gamma n\times \gamma n$ matrix with
eigenvalues $c_{r\gamma+1},\ldots,c_{(r+1)\gamma}$, each of
multiplicity $n$, $Q^j$ are $\R^{\gamma n\times \gamma n}$-valued
matrix polynomials in $y^{(1,\ldots,(j+1)\gamma-1)}$ and $S_{t,r}$
is a $\R^{\gamma n}$-valued polynomial in
$y^{(1,\ldots,(r-t)\gamma)}$ depending on both $t$ and $r$.
Here we put all the terms containing $y^{(i)}$ with $i\le (r-t)\gamma$
into
$S_{t,r}$.

We first try to solve the system (\ref{collect}) with respect to
the $\gamma$ highest derivatives $y^{(r\gamma+1)},\ldots,
y^{((r+1)\gamma)}$. We can do this provided the coefficient matrix
$Q^0 (y^{(1,\ldots,\gamma-1)},\t)$ is invertible. In general, a
solution can be obtained only modulo the kernel of
$Q^0(y^{(1,\ldots,\gamma-1)},\t)$.
Here the dimension of the kernel may change as
$(y^{(1,\ldots,\gamma-1)},\t)$
changes.
To avoid this problem we consider only solutions $y(x,\t)$
of (\ref{eq-m}) with $y^{(1,\ldots,\gamma-1)}= \hat
y^{(1,\ldots,\gamma-1)}$ (as we may by {\it a priori} assuming
$k\geq \gamma-1$).
For these solutions, we obtain from (\ref{collect}) an identity
\begin{equation}\Label{ker}
y^{(r\gamma+1,\ldots,(r+1)\gamma)} = T_{r,t}(y^{(1,\ldots,r\gamma)},\t)
\mod \ker Q^0
\end{equation}
where $T_{r,t}$ is a $\R^{\gamma n}$-valued polynomial in
$y^{(1,\ldots,r\gamma)}$
and $Q^0:=Q^0(\hat y^{(1,\ldots,\gamma-1)},\t)$.
Here two cases are possible.
If the kernel of $Q^0$ is trivial for some $\t$,
it is trivial for $\t$ outside a proper real-analytic subvariety.
Then for such values of $\t$, (\ref{ker}) can be iterated
to determine all derivatives $y^{(s)}$, for $s\geq \gamma+1$, in terms
of $y^{(1,\ldots,\gamma)}$. The proof is complete by continuity.

If the kernel of $Q^0$ is nontrivial for all $\t$,
it has a constant dimension for $\t$ outside a proper
real-analytic subvariety. Then we consider the system
(\ref{collect}) with $r$ replaced by $r+1$. Here the
$\gamma$-tuple of new unknown derivatives
$y^{((r+1)\gamma+1,\ldots,(r+2)\gamma)}$ with the coefficient
matrix $Q^0$ is involved. However, we
can still extract some information when the image $\Im
Q^0$ is a proper subspace of $\bR^{\gamma
n}$ (which happens precisely when $\ker Q^0\ne \{0\}$), namely
\begin{multline}\Label{im1}
(C_{r+1} - Q^1(y^{(1,\ldots,2\gamma-1)},\t)) y^{(r\gamma+1,\ldots,
(r+1)\gamma)} = \\ \sum_{2\le j\le t} Q^j
(y^{(1,\ldots,(j+1)\gamma-1)},\t) y^{((r-j+1)\gamma+1,\ldots,
(r-j+2)\gamma)} +
S_{r+1,t}(y^{(1,\ldots,(r-t+1)\gamma)},\t)
+ Q^0 y^{((r+1)\gamma+1,\ldots,(r+2)\gamma)}.
\end{multline}

We now use the explicit
form of $C_{r+1}$ to conclude that, for every $r$ sufficiently large,
the matrix
\begin{equation}\label{matr}
C_{r+1} - Q^1(y^{(1,\ldots,2\gamma-1)},\t)
\end{equation}
on the left-hand side is invertible for
$y^{(1,\ldots,2\gamma-1)} = \hat y^{(1,\ldots,2\gamma-1)}$
and for $\t$ outside a proper subvariety.
By taking the union of these subvarieties for different $r$
we see that, for each $\t$ outside a countable union of proper
subvarieties, the matrices (\ref{matr}) are invertible for all $r$.
We write
\begin{equation*}
A^1_{r+1}(y^{(1,\ldots,2\gamma-1)},\t):=
\big(C_{r+1} - Q^1(y^{(1,\ldots,2\gamma-1)},\t)\big)^{-1}.
\end{equation*}
Thus $A^1_{r+1}$ is a rational function in $y^{(1,\ldots,2\gamma-1)}$.
By applying $A^1_{r+1}(y^{(1,\ldots,2\gamma-1)},\t)$ to both sides of
(\ref{im1}), we obtain a rational expression for
$y^{(r\gamma+1,\ldots,(r+1)\gamma)}$
in terms of lower order derivatives
modulo the linear subspace
\begin{equation}\Label{Vr}
V^1_{r}(y^{(1,\ldots,2\gamma-1},\t):=
A^1_{r+1}(y^{(1,\ldots,2\gamma-1)},\t) \ \Im Q^0
\subset \R^{\gamma n}.
\end{equation}
Previously we fixed the first derivatives $y^{(1,\ldots,\gamma-1)}$.
In this step, we further assume that
$y^{(1,\ldots,2\gamma-1)}=\hat y^{(1,\ldots,2\gamma-1)}$.
We can then drop the dependence of $Q^1$, $A^1_{r+1}$ and $V^1_r$
on the derivatives as we did for $Q^0$.
Taking also (\ref{ker}) into account, we
conclude that $y^{(r\gamma+1,\ldots, (r+1)\gamma)}$ is determined
modulo
\begin{equation}\Label{intersection}
\ker Q^0\cap V^1_{r}.
\end{equation}
If this space is zero-dimensional (for some $\t$), the two
equations (\ref{ker})
and (\ref{im1}) determine \break
 $y^{(r\gamma+1,\ldots,(r+1)\gamma)}$
completely and yield a polynomial expression for these derivatives
in terms of lower derivatives (for this value of $\t$).

We now observe that we may write
$A^1_{r+1}=c_{(r+2)\gamma}^{-1}B^1_{r+1}$, where the matrix
$B^1_{r+1}$ tends to the identity (for $\t$ fixed)
as $r\to \infty$. Moreover the linear operator
$B^1(\epsilon):=B^1_{r+1}$, where
$\epsilon:=1/r$ (and the integers $r$ in
$B^1_{r+1}$  are replaced by a continuous
variable $r$ in the obvious way),
is a ratio of analytic functions for $(\epsilon,\t)$
in a neighborhood of $0$.
By Cramer's rule,
there exist real-analytic functions
$v_1(\t),\ldots, v_\mu(\t)$ and $u_{\mu+1}(\t),\ldots, u_{n\gamma}(\t)$
that represent bases for $\ker Q^0$
and for $\Im Q^0$ respectively
for every $\t$ outside a proper subvariety.
If we write
$\Delta(\epsilon,\t)$ for the determinant of the matrix of
$n\gamma$-vectors \begin{equation}\Label{det}
v_1,\ldots, v_{\mu},
B^1u_{\mu+1},\ldots, B^1u_{n\gamma},
\end{equation}
then the intersection (\ref{intersection}) is positive dimensional if
and only if $\Delta(\epsilon,\t)$ vanishes at $(1/r,\t)$. Since
$\Delta(\epsilon,\t)$ is real-analytic near $0$, we
conclude that either $\Delta(\epsilon,\t)$ is identically 0 or
there exists $r_0$ such that $\Delta^1(1/r,\t)\neq 0$ for $r\ge r_0$
and for $\t$ outside a proper subvariety $\Gamma_r^1$.
In the second case the intersection (\ref{intersection})
is zero-dimensional for $r\ge r_0$ and $\t\notin \Gamma_r^1$
and we obtain $y^{(r\gamma+1,\ldots,(r+1)\gamma)}$
as a polynomial expression in lower derivatives.
By a simple inductive argument and the analyticity
of $y$ and $\hat y$, it follows that $y(x,\t)=\hat y(x,\t)$ for all $x$
and all $\t$ outside the union of $\Gamma_r^1$, $r\ge r_0$.
The desired statement follows by continuity,
since the union of $\Gamma_r^1$'s is nowhere dense.

Thus, we may assume that $\Delta(\epsilon,\t)\equiv 0$.
In this case, we consider (\ref{collect}) with $r$ replaced by
$r+2$. By solving (\ref{im1}), with $r+1$ in place of $r$,
for $y^{((r+1)\gamma+1,\ldots,(r+2)\gamma)}$
(modulo $\Im Q^0(y^{1,\ldots, \gamma-1},\t)$), we conclude that
\begin{multline}\Label{im2}
y^{((r+1)\gamma+1,\ldots,
(r+2)\gamma)} =
A^1_{r+2}
\Big(
\sum_{2\le j\le t} Q^j
(y^{(1,\ldots,(j+1)\gamma-1)},\t) y^{((r-j+2)\gamma+1,\ldots,
(r-j+3)\gamma)} + \\ S_{r+2,t}(y^{(1,\ldots,(r-t+1)\gamma)},\t)\Big)
+
A^1_{r+2} Q^0 y^{((r+2)\gamma+1,\ldots,(r+3)\gamma)},
\end{multline}
where $A^1_{r+2}$ and
$S_{r+2,t}$ are defined as above.
By applying
$Q^0$ to both sides and substituting
the right-hand side for
$Q^0 y^{((r+1)\gamma+1,\ldots,(r+2)\gamma)}$
in (\ref{im1}), we deduce that
\begin{equation}\Label{im3}
D^2_{r+1}(y^{(1,\ldots,3\gamma-1)},\t)y^{(r\gamma+1,\ldots,
(r+1)\gamma)} =
 R^2_{r+1}(y^{(1,\ldots,r\gamma)},\t) \mod
Q^0 (V^1_{r+1}),
\end{equation}
where $V^1_{r+1}$ is as defined above and $D^2_{r+1}$ is the
invertible (for $r$
large enough and $\t$ outside a proper subvariety $\Gamma_r^2$,
provided that we have $y^{(1,\ldots,3\gamma-1)}=\hat
y^{(1,\ldots,3\gamma-1)}$) matrix
\begin{equation}\Label{Dr}
D^2_{r+1}(y^{(1,\ldots,3\gamma-1)}):=C_{r+1} -
Q^1
-Q^0
A^1_{r+2}
Q^2
(y^{(1,\ldots,3\gamma-1)},\t).
\end{equation}
As before we assume in this step
that $y^{(1,\ldots,3\gamma-1)}= \hat y^{(1,\ldots,3\gamma-1)}$
and drop the dependence on these derivatives.

Now observe that the assumption that
$\Delta(\epsilon,\t)\equiv 0$ (which is equivalent to the
intersections (\ref{intersection}) being nontrivial for all $r$
sufficiently large) implies that the space
$Q^0 (V^1_{r+1})$ is of
strictly lower dimension than $V^1_{r+1}$.
This is a crucial observation. Let us
write $A^2_{r+1}$ for the inverse of
$D^2_{r+1}$. It follows that we can
solve (\ref{im3}) for $y^{(r\gamma+1,\ldots, (r+1)\gamma)}$ modulo
the subspace
\begin{equation}\Label{V2r}
V^2_r
:=A^2_{r+1}
Q^0
(V^1_{r+1}).
\end{equation}
If the
intersection
\begin{equation}\Label{intersection2}
\ker Q^0
\cap
V^2_{r}
\end{equation}
is zero-dimensional, we can find a polynomial expression for
$y^{(r\gamma+1,\ldots,(r+1)\gamma)}$  in terms of lower
derivatives by using equation (\ref{ker}). We claim again that
this intersection will be either zero-dimensional for all
sufficiently large $r$ and $\t$ outside a proper subvariety
$\Gamma^3_r$
or positive-dimensional for all
sufficiently large $r$ and all $\t$. The argument is as before.  We can
detect
positive dimensionality of the intersection (\ref{intersection2})
by the vanishing of suitable determinants formed by the vectors in
(\ref{det}) with $B^1$ replaced by $B^2=B^2(\epsilon)$,
where $B^2(\epsilon)$ (which also of course depends on $\t$) is defined
as
follows. Let us factor the
scalar $c_{(r+2)\gamma}^{-1}$ in
$A^2_{r+1}$, writing
$A^2_{r+1}
=c_{(r+2)\gamma}^{-1}
B^2_{r+1}$.
Then we define
$B^2(\epsilon):= B^2_{r+1} Q^0 B^1_{r+2}$,
where as before $\epsilon:=1/r$. It
is not difficult to see that $B^2$ is analytic in $(\epsilon,\t)$ near
$0$
with $B^2(0,\t)$ equal to the identity. Since any
determinant formed by the vectors in (\ref{det}) with
$B^1$ replaced by $B^2$ will be analytic,
the claim now follows as above.

As mentioned above, if
the intersection (\ref{intersection2}) is trivial for all sufficiently
large $r$, we are done.
If not, we must go iterate the procedure
above, and start with the equation (\ref{collect}) with $r$
replaced by $r+3$. In this way we will obtain a subspace
$V^3_r(y^{(1,\ldots,4\gamma-1)},\t)$ (in a way analogous to that
yielding
$V^2_r(y^{(1,\ldots,3\gamma-1)},\t)$). By the same argument as above,
the fact that we are forced to go to the next iteration (i.e.\ the
intersection (\ref{intersection2}) is nontrivial for all large
$r$) implies that $V^3_r$ has strictly
lower dimension than $V^2_r$. The
crucial observation is that, if we are forced to make another
iteration, the dimension of the subspaces
$V^j_r:=V^j_r(\hat y^{(1,\ldots,(j+1)\gamma-1)},\t)$ drops. Hence, the
process
will terminate after at most $n\gamma$ steps. The details of the
iterations are left to the reader.

Summarizing, we obtain a linear
system for
$y^{(r\gamma+1,\ldots,(r+1)\gamma)}$ (providing $r$ is large enough)
in terms of lower order derivatives
and the matrix coefficient of $y^{(r\gamma+1,\ldots,(r+1)\gamma)}$ is
polynomial
in $y^{(1,\ldots,(n\gamma+1)\gamma)}$
and is invertible for $\t$ outside a countable union of proper
subvarieties.
Then the proof is completed by the analyticity of $y(x,\t)$ and $\hat
y(x,\t)$
and by continuity as before.
\end{proof}

\section{Finite jet determination in the infinite type case;
proof of Theorem~\ref{main} }\Label{sec-infinite}
In this section we complete the proof of Theorem~\ref{main}.
If $M$ is of finite type at $p$, the statement is a special case of
Theorem~\ref{main-jetpar}.
Hence, to complete the proof of Theorem \ref{main}, we may assume
that the hypersurface $M\subset\C^2$
is of {\em infinite type} at $p$.
This is equivalent to the property that the Segre variety $E$ of $p$
is contained in $M$.
As before we denote the same objects
associated to another real-analytic hypersurface $M'\subset\C^2$ by
$'$.
Given $M'$, $H^1$ and $H^2$ as in Theorem~\ref{main},
we set $p':=H^1(p)$. It follows that $M'$ is also of infinite type at
$p'$.
The proof of Theorem \ref{main} in this case is based on
Theorem~\ref{thm-jetsegre}
and on the following result.

\begin{Thm}\Label{E-cor}
Let $h^0$ be a $C^\infty$-smooth CR-diffeomorphism between
real-analytic hypersurfaces $M$ and $M'$ in $\C^2$. Suppose that
$M$ is of infinite type at a point $p\in M$ and set $p':=h^0(p)\in
M'$. Choose local coordinates $y=(x,s)\in\bR^{2}\times \bR$ on $M$
and $y'=(x',s')\in \bR^{2}\times\bR$ on $M'$ vanishing at $p$ and
$p'$, respectively, such that the zeroth order Segre varieties
$E\subset M$ and $E'\subset M'$ at $p$ and $p'$ are locally given
by $s=0$ and $s'=0$, respectively. Then there exists an integer
$m\ge 1$ such that, if $h$ is a $C^\infty$-smooth
CR-diffeomorphism between open neighborhoods of $p$ and $p'$ in
$M$ and $M'$ respectively with $h(p)$ sufficiently close to $p'$
and we set $g(y):=s'(h(y))$, then there is a (unique)
$C^\infty$-smooth function $v_h(y)$ on $M$ near $p$ satisfying
\begin{equation}\Label{system-0}
s^m\partial_s g(y) \equiv v_h(y) g(y)^{m}.
\end{equation}
For any such $h$ we write
$$u_h(y):=\big((\partial_{x_i} h(y))_{1\le i\le 2}, s^m\partial_s f(y),
v_h(y))\in \R^9,$$
where $f(y):=x'(h(y))$, and set
\begin{equation} \Label{eq-contact}
u^{\a}_h(y):=(s^m\partial_s)^{\a_0}\d_{x_1}^{\a_1}\d_{x_2}^{\a_2} u_h
(y)
\end{equation}
for any multi-index $\a=(\a_0,\a_1,\a_2)\in \mathbb Z_+^3$.
Then there exist real-analytic
functions $q(y)$ on $M$ near $p$
and $q'(y')$ on $M'$ near $p'$ with $q(x,0)$ and $q'(x',0)$ both not
identically zero, an open neighborhood
$\Omega\subset J^2(M,\bR^9)_{p}\times M'\times M$
of $\big((u^\b_{h^0}(p))_{0\le|\b|\le 2},p',p\big)$,
and, for every multi-index $\a\in \mathbb Z_+^3$ with
$|\a|=3$,
real-analytic functions $r^\a(\L,y',y)$ on $\Omega$ such that,
for any $h$ as above with
\begin{equation}\Label{h-cond}
\big((u^\b_{h}(p))_{0\le|\b|\le 2},h(p),p\big)\in\Omega,
\end{equation} the equation
\begin{equation}\Label{system}
u_h^\a(y) \equiv
\frac{r^{\a}\big((u_h^\b(y))_{0\le|\b|\le2},h(y),y\big)}
{q(y) q'(h(y))}
\end{equation}
holds at every $y\in M$ near $p$ for which  the denominator does not
vanish.
\end{Thm}

At points $p\in M$ and $p'\in M'$ where $q(p)\neq 0$ and
$q'(p')\neq 0$ (the latter of which is, in fact, a consequence of
the former), Theorem \ref{E-cor} is a reformulation of Theorem 2.1
in \cite{Enonm}. Theorem \ref{E-cor}, as stated above, follows by
repeating the proof of Theorem 2.1 in \cite{Enonm} at a general
point p, accepting the presence of a denominator which may vanish
at $p$. The result is the conclusion of Theorem \ref{E-cor}. We
shall omit the details and refer the reader to the proof in
\cite{Enonm} for inspection.

\bpf[Proof of Theorem~\ref{main}] Let $H^1,H^2$ be as in
Theorem~\ref{main} with $j^k_p H^1=j^k_pH^2$, for some $k\geq 2$
(which will be specified later).  Observe that, for any local
biholomorphism $H\colon (\C^2,p)\to \C^2$ sending $M$ into itself,
the restriction $h:=H|_M$ is a local CR-diffeomorphism between
open neighborhoods of $p$ and $p':=H(p)$ in $M$. We shall write
$h^j:=H^j|_M$, $j=1,2$. If we take $h^0:=h^1$, $h:=h^2$ and $M'=M$
in Theorem \ref{E-cor}, then $h$ satisfies (\ref{h-cond}) and
hence the equation (\ref{system}). For an $h$ as in
Theorem~\ref{E-cor}, we set
$$U_h(y):=\big((u^\b_{h}(y))_{0\le|\b|\le 2},h(y)\big).$$
Equations  (\ref{system-0}) and (\ref{system}) then imply that
\begin{equation}\label{R-id}
s^m\d_s U_h(y) \equiv \frac{R(U_h,y)}{q(y) q'(h(y))}
\end{equation}
for some real-analytic function $R(U,y)$
defined in a neighborhood of $(U_{h^0},0)$.

The assumption $j^k_p H^1 \equiv j^k_p H^2$
and Theorem~\ref{thm-jetsegre} imply that
\begin{equation}\Label{eq-jetatsegre}
j^{k-1}_Z H^1 \equiv j^{k-1}_Z
H^2,\quad\forall Z\in E,\end{equation}
near $p$.
Let us write, for $j=1,2$,
 \begin{equation}\Label{eq-expU}
U_{h^j}(x,s)=\sum_{k=0}^\infty U^j_k(x)s^k.
\end{equation}
We conclude, by the construction of $U_h$ and (\ref{eq-jetatsegre}),
that $U^1_{l}(x) \equiv U^2_{l}(x)$ for $l\leq k-4$. .
To complete the proof of Theorem \ref{main}, consider the
functions $\tilde U_{h^j}(x,s):=U_{h^j}(x,s)-U_{h^j}(x,0)$, for
$j=1,2$,
which satisfy $\tilde U_{h^j}(x,0)=0$.
For $k\ge 4$, we have $U_{h^1}(x,0)=U_{h^2}(x,0)$ and hence both
$\tilde U_{h^1}$ and $\tilde U_{h^2}$ satisfy the same system of
differential
equations
\begin{equation}\Label{eq-system2}
s^m\partial_s
  \tilde U_h(y)=\frac{\tilde R(y,\tilde U_h)}{q(y) q'(h(y))},
\end{equation}
where $\tilde R(y,\tilde U):=R(y,\tilde U+U^1(x,0))$, as does any
other $\tilde U_h$ arising from a CR diffeomorphism $h$ with
$U_h(x,0)=U^1(x,0)$. Recall that $h(y)$ is one of the components
of $U_h$. The existence of the integer $k$ such that $H^1\equiv
H^2$ (which is equivalent to $h^1\equiv h^2$) if
$j^{k}_{p}H^1=j^{k}_{p}H^2$ now follows from Theorem~\ref{ode},
although the choice of the integer $k$ appears to depend on the
mapping $H^1$. However, we shall show that one can find a $k$ that
works for every $H^1$. Let $k$ be the integer obtained by the
  above procedure applied to $H^1,H^2$ where
  $H^1(Z)$ is the identity mapping $\id$.
We conclude that if $H\colon (\C^2,p)\to (\C^2,p)$
sends $M$ into itself and $j^{k}_p H=\id$, then
  $H\equiv \id$. We claim that the same number $k$
satisfies the conclusion of Theorem~\ref{main} for any  $H^1$,
$H^2$. Indeed, for any $H^1$, $H^2$  as in Theorem~\ref{main}, the
mapping $H:=(H^1)^{-1}\circ H^2$ sends $(M,p)$ into itself and
satisfies $j^{k}_p H=\id$. Hence, by
  the construction of $k$, we must have $(H^1)^{-1}\circ H^2\equiv \id$
which
  proves the claim. This completes the proof of Theorem~\ref{main}.
\epf


\begin{thebibliography}{CNWS99}



\bibitem[BER97]{Asian} Baouendi,~M.S.; Ebenfelt,~P.; Rothschild,~L.P.:
Parametrization of local biholomorphisms of
real-analytic hypersurfaces, {\em Asian J. Math.} {\bf 1}, 1--16,
(1997).


\bibitem [BER98] {CAG} Baouendi, M. S.; Ebenfelt, P.; Rothschild, L.P.:
CR automorphisms of real analytic manifolds in
complex space. {\em Comm. Anal. Geom.} {\bf 6},  291--315, (1998).


\bibitem[BER99a]{BER} Baouendi,~M.S.; Ebenfelt,~P.; Rothschild,~L.P.:
{\em Real Submanifolds in Complex Space and Their Mappings}.
Princeton Math. Series {\bf 47}, Princeton Univ. Press, 1999.


\bibitem[BER99b]{MA} Baouendi,~M.S.; Ebenfelt,~P.; Rothschild,~L.P.:
Rational dependence of smooth
and analytic CR mappings on their jets. {\em Math.~Ann.}
{\bf 315}, 205--249, (1999).


\bibitem [BER00a]{JAMS} Baouendi,~M.S.; Ebenfelt,~P.; Rothschild,~L.P.:
Convergence and finite determination of formal CR mappings.
{\em J. Amer. Math. Soc.} {\bf 13}, 697--723, (2000).


\bibitem[BMR00]{BMR} Baouendi, M.S.; Mir, N.; and Rothschild, L.P.:
Reflection ideals and mappings between generic submanifolds in complex
space. (preprint 2000),
{\tt <http://xxx.lanl.gov/abs/math.CV/0012243>}.


\bibitem[Be90]{Be90} Beloshapka, V.K.:
A uniqueness theorem for automorphisms of a nondegenerate surface in a
complex space.
{\em Math. Notes} {\bf 47}, No.3, 239--242 (1990);
translation from {\em Mat. Zametki} {\bf 47}, No.3, 17--22, (1990).


\bibitem[BG77]{BG77} Bloom,~T.; and Graham,~I.:
On type conditions for generic real submanifolds of {$\bC^n$}. {\em
Invent.~Math.} {\bf 40}, 217--243, (1977).


\bibitem[CaE32a]{Ca32a} Cartan,~E.: Sur la g\'eom\'etrie
pseudo-conforme des
hypersurfaces de deux variables complexes, I. {\em Ann. Math. Pura
Appl.} {\bf
11}, 17--90, (1932). ({\em \OE uvres compl\`etes}, Part. II, Vol. 2,
Gauthier-Villars, 1231--1304, (1952)).


\bibitem[CaE32b]{Ca32b} Cartan,~E.: Sur la g\'eom\'etrie
pseudo-conforme des
hypersurfaces de deux variables complexes, II. {\em Ann. Sc. Norm. Sup.
Pisa}
{\bf 1}, 333--354, (1932). ({\em \OE uvres compl\`etes}, Part. III,
Vol. 2,
Gauthier-Villars, 1217--1238, (1952)).


\bibitem[CaH35]{CaH35} Cartan,~H.:
{\em Sur les groupes de transformations analytiques.} Act.~Sc.~et~Int.,
Hermann,~Paris, 1935.


\bibitem[CM74]{CM} Chern, S.-S.; Moser, J.K.:
Real hypersurfaces in complex manifolds. {\em Acta Math. } {\bf 133},
219--271, (1974).


\bibitem[D82]{D82} D'Angelo,~J.P.: Real hypersurfaces, orders of
contact, and applications.
{\em Ann.~Math.} {\bf 115}, 615--637, (1982).


\bibitem[E00a]{Ejet} Ebenfelt, P.:
Finite jet determination of holomorphic mappings at the boundary.
{\em Asian J. Math.} (to appear, 2000),
{\tt <http://xxx.lanl.gov/abs/math.CV/0001116>}.


\bibitem[E00b]{Enonm} Ebenfelt, P.:
On the analyticity of CR mappings between nonminimal hypersurfaces.
(preprint 2000),
{\tt <http://xxx.lanl.gov/abs/math.CV/0009010>}.


\bibitem [Han97]{Han2} Han, C.K.: Complete system for the mappings
of CR manifolds of nondegenerate Levi forms.
{\em Math. Ann.} {\bf 309}, 401--409, (1997).


\bibitem [Hay98]{Hay} Hayashimoto, A.: On the complete system of
finite order for CR mappings and its application. {\em Osaka J. Math.}
{\bf 35},  617--628, (1998).


\bibitem [Hu96]{HCPDE} Huang, X.: Schwarz reflection
principle in complex spaces of dimension two. {\em
Comm. Partial Diff. Eq.} {\bf  21},  1781--1828, (1996).


\bibitem[J77]{J77} Jacobowitz,~H.:
Induced connections on hypersurfaces in ${\C}\sp{n+1}$.
{\em Invent. Math.} {\bf 43}, no. 2, 109--123, (1977).


\bibitem[Ki01]{Ki01} Kim,~S.-Y.:
Complete system of finite order for CR mappings between
real analytic hypersurfaces of degenerate Levi form.
{\em J. Korean Math. Soc.} {\bf 38}, no. 1, 87--99, (2001).


\bibitem[KZ01]{KZ} Kim,~S.-Y.; and Zaitsev,~D.:
The equivalence problem for CR-structures of any codimension.
(preprint 2001).


\bibitem[Ko72]{Ko72} Kohn,~J.J.:
Boundary behavior of $\bar \partial$ on weakly pseudo-convex
manifolds of dimension two. {\em J. Differential Geom.} {\bf 6},
523--542, (1972).


\bibitem[L01]{Lamel} Lamel, B.: Holomorphic
maps of real submanifolds in complex spaces of
different dimensions. {\em Pacific J. Math.} (to appear, 2001),
{\tt <http://xxx.lanl.gov/abs/math.CV/9911057>}.


\bibitem[Ta62]{Ta62} Tanaka,~N.: On the pseudo-conformal geometry of
hypersurfaces of the space of $n$ complex variables. {\em J. Math. Soc.
Japan}
{\bf 14}, 397--429, (1962).


\bibitem[To66]{To66} Tomassini,~G.:
Tracce delle functional olomorfe sulle sotto varieta analitiche reali
d'una varieta complessa.
{\em Ann.~Scuola~Norm.~Sup.~Pisa} {\bf 20} (1966), 31--43.


\bibitem[Z97]{Z1} Zaitsev, D.: Germs of
local automorphisms of real-analytic CR structures and analytic
dependence on $k$-jets.  {\em Math. Res. Lett.} {\bf 4},
823--842, (1997).


\end{thebibliography}
\end{document}